\documentclass[12pt]{article}
\usepackage[english,russian]{babel}
\usepackage{amsfonts,amssymb,amscd,amsmath,amsthm}
\def\R{{\bf R}}   
\long\def\comment#1\endcomment{}

\def\id{\mathop{\fam0 id}}

\def\im{\mathop{\fam0 im}}

\def\eps{\varepsilon}
\usepackage{graphicx}
\def\mpfile#1#2{\includegraphics{#1#2.eps}}

\begin{document}

\newpage
\centerline{\uppercase{\bf ┴рчшёэ√х тыюцхэш  ш 13-  яЁюсыхьр ├шы№схЁЄр}
\footnote{
╬Єфхы№э√х ўрёЄш ЁрсюЄ√ яЁхфёЄрты ышё№ эр ыхЄэшї ъюэЇхЁхэЎш ї ╥єЁэшЁр
├юЁюфют 1997 ш 2006 уюфют (┬.└. ├юЁшэ√ь, ┬.└. ╩єЁышэ√ь, ╚.═. ╪эєЁэшъют√ь
ш ртЄюЁюь), эр ёхьшэрЁрї ьхїьрЄр ╠├╙ ш эр ёхьшэрЁх яю ухюьхЄЁшш т ╠╓═╠╬.
═хъюЄюЁ√х чрфрўш ю уырфъющ срчшёэюёЄш яЁхфёЄрты ышё№ └.└. ┴рЁрэюь эр
ьхцфєэрЁюфэющ ъюэЇхЁхэЎшш °ъюы№эшъют Intel ISEF т 2003 уюфє.
╤ё√ыъш фр■Єё  яю тючьюцэюёЄш эх эр юЁшушэры№э√х ЁрсюЄ√, р эр юсчюЁ√.}
}
\smallskip
\centerline{\bf └. ╤ъюяхэъют}
\smallskip


\bigskip

┬ ¤Єющ ёЄрЄ№х Ёрёёърчрэю, ъръ яЁш Ёх°хэшш 13-щ яЁюсыхь√ ├шы№схЁЄр ю
ёєяхЁяючшЎш ї эхяЁхЁ√тэ√ї ЇєэъЎшщ яю тшыюё№ яюэ Єшх срчшёэюую яюфьэюцхёЄтр ш
срчшёэюую тыюцхэш .
╧Ёштюф Єё  ъЁрёшт√х Ёхчєы№ЄрЄ√ юс ¤Єшї яюэ Єш ї, сюы№°р  ўрёЄ№ ъюЄюЁ√ї
фюёЄєяэр ёЄрЁ°хъырёёэшъє.
╥Ёш ўрёЄш ёЄрЄ№ш ьюцэю ўшЄрЄ№ эхчртшёшью фЁєу юЄ фЁєур (т эхсюы№°юь ъюышўхёЄтх
ьхёЄ, т ъюЄюЁ√ї юфэр ўрёЄ№ шёяюы№чєхЄ фЁєує■, яЁштхфхэ√ Єюўэ√х ёё√ыъш).

┬ яхЁтющ ўрёЄш яЁштюфшЄё  ¤ыхьхэЄрЁэюх шчыюцхэшх шфхш Ёх°хэш  13-щ яЁюсыхь√
├шы№схЁЄр ю ёєяхЁяючшЎш ї └. ═. ╩юыьюуюЁют√ь ш ┬. ╚. └Ёэюы№фюь (яю ьюЄштрь
[Ar58]).
╧Ёш ¤Єюь яюърч√трхЄё , ъръ хёЄхёЄтхээю яю тшыюё№ яюэ Єшх срчшёэюую тыюцхэш , ш
юёЄрхЄё  схч фюърчрЄхы№ёЄтр трцэхщ°хх ьхёЄю (ыхььр ╩юыьюуюЁютр ю фхЁхт№ ї).

╧юфьэюцхёЄтю $K\subset\R^2$ эрч√трхЄё  {\it срчшёэ√ь}, хёыш
фы  ы■сющ эхяЁхЁ√тэющ ЇєэъЎшш $f:K\to\R$ ёє∙хёЄтє■Є Єръшх эхяЁхЁ√тэ√х ЇєэъЎшш
$g,h:\R\to\R$, ўЄю $f(x,y)=g(x)+h(y)$ фы  ърцфющ Єюўъш $(x,y)\in K$.

┬ю тЄюЁющ ўрёЄш яЁштюфшЄё  їрЁръЄхЁшчрЎш  уЁрЇют,
ъюЄюЁ√х ьюцэю тыюцшЄ т яыюёъюёЄ№ т ърўхёЄтх срчшёэ√ї яюфьэюцхёЄт
[St89, Sk95] (Ёх°хэшх яЁюсыхь√ ╪ЄхЁэЇхы№фр),
р Єръцх хх юсюс∙хэш  [St89, Sk95, Ku00, Ku03, Ku03'].

╥ЁхЄ№  ўрёЄ№ эршсюыхх ¤ыхьхэЄрЁэр (ёь., эряЁшьхЁ, чрфрўш 1b ш 4a).
╧ЁштюфшЄё  їрЁръЄхЁшчрЎш  срчшёэ√ї яюфьэюцхёЄт яыюёъюёЄш (Ёх°хэшх яЁюсыхь√
└Ёэюы№фр) [Ar58', St89, MKT03, Mi09], р Єръцх хх юсюс∙хэш  [St89, Vo81, Vo82,
RZ07].

┬ ЄхъёЄх ьэюую чрфрў, ъюЄюЁ√х эєьхЁє■Єё  цшЁэ√ьш ЎшЇЁрьш ш ъ сюы№°шэёЄтє
ъюЄюЁ√ї яЁштюф Єё  Ёх°хэш  (т ъюэЎх ёююЄтхЄёЄтє■∙хую яєэъЄр).
┼ёыш єёыютшх чрфрўш  ты хЄё  єЄтхЁцфхэшхь, Єю чрфрўр ёюёЄюшЄ т Єюь, ўЄюс√
¤Єю єЄтхЁцфхэшх фюърчрЄ№.
╥Ёєфэ√х чрфрўш юЄьхўхэ√ чтхчфюўъющ, р эхЁх°хээ√х --- фтєь .

┼ёыш эхъюЄюЁюх чрьхўрэшх т ёэюёъх шыш єёыютшх чрфрўш ┬рь эхяюэ Єэю, Єю хую
эєцэю яЁюёЄю шуэюЁшЁютршЄ№.
▌Єю эх яЁштхфхЄ ъ эхяюэшьрэш■ фры№эхщ°хую ьрЄхЁшрыр.

┴ыруюфрЁ■ ┬.╚. └Ёэюы№фр, ▐.╠. ┴єЁьрэр, ╤.╠. ┬юЁюэшэр, ╠.═. ┬ ыюую, └.╨. ╤рЇшэр 
ш ╚.═. ╪эєЁэшъютр чр яюыхчэ√х юсёєцфхэш , 
р Єръцх ╠.═. ┬ ыюую чр яюфуюЄютъє Ёшёєэъют.


\bigskip
\centerline{\uppercase{\bf ╬ Ёх°хэшш 13-щ яЁюсыхь√ ├шы№схЁЄр}}

\smallskip
{\bf 13-  яЁюсыхьр ├шы№схЁЄр.}

╧єёЄ№ фрэю эхъюЄюЁюх ьэюцхёЄтю ЇєэъЎшщ
$F= \{f_\alpha(x_1,\dots,x_{n_\alpha})\}_{\alpha\in A}$
(эх юс чрЄхы№эю ъюэхўэюх).
╬яЁхфхышь {\it ёєяхЁяючшЎш■} ЇєэъЎшщ шч $F$ (шыш {\it ЇюЁьєыє} эрф $F$)
шэфєъЄштэю:

(1) ёрьш ЇєэъЎшш $f_\alpha$ ш тёх яхЁхьхээ√х $x_j$  ты ■Єё  ёєяхЁяючшЎш ьш
ЇєэъЎшщ шч $F$.

(2) хёыш ЇєэъЎшш $f(x_1,\dots,x_n),\ g_1(\dots),\ \dots,\
g_n(\dots)$  ты ■Єё  ёєяхЁяючшЎш ьш ЇєэъЎшщ шч $F$ (эх юс чрЄхы№эю
Ёрчышўэ√ьш), Єю ш ЇєэъЎш  $f(g_1(\dots),\dots,g_n(\dots))$
 ты хЄё  ёєяхЁяючшЎшхщ ЇєэъЎшщ шч $F$.
╟фхё№ т ърўхёЄтх рЁуєьхэЄют ЇєэъЎшщ $g_i$ ьюцэю яюфёЄрты Є№ ы■с√х эрсюЁ√ 
яхЁхьхээ√ї (¤Єш яхЁхьхээ√х ьюуєЄ шфЄш т ы■сюь яюЁ фъх, р эхъюЄюЁ√х шч 
яхЁхьхээ√ї ьюуєЄ ёютярфрЄ№). 
\footnote{╬яЁхфхыхэшх ёєяхЁяючшЎшш ьюцэю Єръцх ёЇюЁьєышЁютрЄ№ уЁрЇшўхёъш,
эр  ч√ъх ёїхь.}

═ряЁшьхЁ, яюышэюь $\sum a_{k_1\ldots k_n}x_1^{k_1}\cdots x_n^{k_n}$
 хёЄ№ ёєяхЁяючшЎш  ЇєэъЎшщ $f(x,y)=x+y$ ш $g(x,y)=xy$ ш ъюэёЄрэЄ.
╧Ёшўхь хёыш ьюцэю шёяюы№чютрЄ№ ЇєэъЎшш юфэющ яхЁхьхээющ, Єю фы  $x,y>0$
яЁюшчтхфхэшх ьюцэю ш эх шёяюы№чютрЄ№, Єръ ъръ $xy=2^{\log_2 x+\log_2 y}$.

╨рёёьюЄЁшь ёыхфє■∙шх тюяЁюё√.

{\it 1. ╠юцэю ыш ърцфє■ ЇєэъЎш■ эхёъюы№ъшї рЁуєьхэЄют чряшёрЄ№ т
тшфх ёєяхЁяючшЎшш ЇєэъЎшщ эх сюыхх ўхь фтєї рЁуєьхэЄют?

2. ╠юцэю ыш ърцфє■ ЇєэъЎш■ фтєї рЁуєьхэЄют чряшёрЄ№ ъръ ёєяхЁяючшЎш■ ЇєэъЎшщ
юфэюую рЁуєьхэЄр ш яЁюёЄхщ°хщ ЇєэъЎшш фтєї рЁуєьхэЄют (эряЁшьхЁ, ёыюцхэш )?
\footnote{─ы  ЇєэъЎшщ рыухсЁ√ ыюушъш юЄтхЄ√ эр юср тюяЁюёр яюыюцшЄхы№э√
(яюёъюы№ъє ы■сє■ ЇєэъЎш■ рыухсЁ√ ыюушъш ьюцэю т√ЁрчшЄ№ ўхЁхч `ш' ш `эх').
└яяЁюъёшьрЎшюээр  ЄхюЁхьр ┬хщхЁ°ЄЁрёр яюърч√трхЄ, ўЄю ЇєэъЎш  эхёъюы№ъшї
рЁуєьхэЄют ьюцхЄ с√Є№ {\it ЁртэюьхЁэю яЁшсышцхэр} эр ъюьяръЄэюь ьэюцхёЄтх
яюышэюьрьш, Є.х. ёєяхЁяючшЎш ьш ъюэёЄрэЄ, ёыюцхэш  ш єьэюцхэш .}
}

╧юёъюы№ъє яыюёъюёЄ№ ш яЁ ьр  Ёртэюью∙э√, Єю ы■сє■ ЇєэъЎш■ ЄЁхї ш сюыхх
яхЁхьхээ√ї ьюцэю т√ЁрчшЄ№ т тшфх ёєяхЁяючшЎшш (тююс∙х уютюЁ , {\it ЁрчЁ√тэ√ї})
ЇєэъЎшщ фтєї яхЁхьхээ√ї (ёь. фхЄрыш т [Ar58]).
╧ю¤Єюьє єърчрээ√х тюяЁюё√ шэЄхЁхёэ√ Єюы№ъю фы  {\it эхяЁхЁ√тэ√ї} ЇєэъЎшщ.

╫хЁхч
$$|x,y|=|(x_1,\dots,x_n),\ (y_1,\dots,y_n)|=
\sqrt{(x_1-y_1)^2+\dots+(x_n-y_n)^2}$$ юсючэрўрхЄё 
юс√ўэюх ЁрёёЄю эшх ьхцфє Єюўърьш $x=(x_1,\dots,x_n)$ ш $(y_1,\dots,y_n)$
яЁюёЄЁрэёЄтр $\R^n$.
╧єёЄ№ $K$ --- яюфьэюцхёЄтю яЁюёЄЁрэёЄтр $\R^n$.
╘єэъЎш  $f:K\to\R$ эрч√трхЄё  {\it эхяЁхЁ√тэющ}, хёыш фы  ы■с√ї Єюўъш
$x_0\in K$ ш ўшёыр $\eps>0$ ёє∙хёЄтєхЄ Єръюх ўшёыю $\delta>0$, ўЄю фы  ы■сющ
Єюўъш $x\in K$ ё єёыютшхь $|x,x_0|<\delta$ т√яюыэхэю $|f(x)-f(x_0)|<\eps$.
═ряЁшьхЁ, ЇєэъЎш  $f(x_1,x_2)=\sqrt{x_1^2+x_2^2}$  ты хЄё  эхяЁхЁ√тэющ эр
яыюёъюёЄш, р ЇєэъЎш   $f(x_1,x_2)$, Ёртэр  Ўхыющ ўрёЄш юЄ $x_1+x_2$, --- эхЄ.
┬ фры№эхщ°хь тёх ЇєэъЎшш яЁхфяюырур■Єё  эхяЁхЁ√тэ√ьш, хёыш  тэю эх юуютюЁхэю
яЁюЄштэюх.

▀ёэю, ўЄю ы■ср  ¤ыхьхэЄрЁэр  ЇєэъЎш  яЁхфёЄрты хЄё  т тшфх ёєяхЁяючшЎшш
ЇєэъЎшщ фтєї яхЁхьхээ√ї.
╧ЁюёЄхщ°шх эх¤ыхьхэЄрЁэ√х ЇєэъЎшш --- ъюЁэш рыухсЁршўхёъшї єЁртэхэшщ.
╩ 1900 уюфє с√ыю шчтхёЄэю, ўЄю ы■сюх рыухсЁршўхёъюх єЁртэхэшх $n$-щ ёЄхяхэш 
ётюфшЄё  (яЁш яюью∙ш Ёрфшърыют, ёыюцхэш , т√ўшЄрэш , єьэюцхэш  ш фхыхэш ) 
ъ Єръюьє, є ъюЄюЁюую ъю¤ЇЇшЎшхэЄ√ яЁш $x^n$ ш $x^0$
Ёртэ√ 1, р яЁш $x^{n-1}$, $x^{n-2}$ ш $x^{n-3}$ Ёртэ√ 0. 
╥ръшь юсЁрчюь, юёЄрхЄё  $n-4$ яхЁхьхээ√ї ъю¤ЇЇшЎшхэЄр.
╧ю¤Єюьє `яЁюёЄхщ°хщ' ЇєэъЎшхщ, фы  ъюЄюЁющ эх с√ыю шчтхёЄэю
т√Ёрцхэшх ўхЁхч ЇєэъЎшш фтєї яхЁхьхээ√ї, с√ыр ЇєэъЎш  $x(a,b,c)$,
т√Ёрцр■∙р  Ёх°хэшх єЁртэхэш  $x^7+ax^3+bx^2+cx+1=0$ ёхф№ьющ ёЄхяхэш.
╧ю¤Єюьє ├шы№схЁЄ ёЇюЁьєышЁютры ётю■ 13-■ яЁюсыхьє Єръ:

{\it ─юърчрЄ№, ўЄю єЁртэхэшх ёхф№ьющ ёЄхяхэш $x^7+ax^3+bx^2+cx+1=0$
эхЁрчЁх°шью схч шёяюы№чютрэш  ЇєэъЎшщ ЄЁхї яхЁхьхээ√ї.}

├шы№схЁЄє єфрыюё№ яюърчрЄ№, ўЄю эхъюЄюЁ√х {\it рэрышЄшўхёъшх} ЇєэъЎшш ЄЁхї
яхЁхьхээ√ї эх  ты ■Єё  ёєяхЁяючшЎш ьш рэрышЄшўхёъшї цх ЇєэъЎшщ фтєї
яхЁхьхээ√ї [Ar58].
┬ 1954 ┬шЄє°ъшэ фюърчры, ўЄю эхъюЄюЁ√х {\it $r$ Ёрч эхяЁхЁ√тэю
фшЇЇхЁхэЎшЁєхь√х} ЇєэъЎшш эх  ты ■Єё  ёєяхЁяючшЎш ьш $r$ Ёрч эхяЁхЁ√тэю
фшЇЇхЁхэЎшЁєхь√ї ЇєэъЎшщ фтєї яхЁхьхээ√ї [Ar58, Vi04].
─ы  {\it эхяЁхЁ√тэ√ї} цх ЇєэъЎшщ ушяюЄхчр ├шы№схЁЄр
с√ыр юяЁютхЁуэєЄр т 1957 уюфє ╩юыьюуюЁют√ь ш └Ёэюы№фюь.

\smallskip
{\bf ╥хюЁхьр ╩юыьюуюЁютр-└Ёэюы№фр.} {\it ╦■ср  эхяЁхЁ√тэр  ЇєэъЎш 
яЁхфёЄрты хЄё  т тшфх ёєяхЁяючшЎшш эхяЁхЁ√тэ√ї ЇєэъЎшщ юфэюую ш фтєї
рЁуєьхэЄют.}


\bigskip
{\bf ╥хюЁхьр ╩юыьюуюЁютр: ъ ёєяхЁяючшЎш ь ЇєэъЎшщ ЄЁхї яхЁхьхээ√ї.}

╤эрўрыр т 1956 у. ╩юыьюуюЁютє єфрыюё№ фюърчрЄ№, ўЄю {\it яЁюшчтюы№эр  эхяЁхЁ√тэр 
ЇєэъЎш  сюыхх ўхь ЄЁхї яхЁхьхээ√ї  ты хЄё  ёєяхЁяючшЎшхщ эхяЁхЁ√тэ√ї ЇєэъЎшщ
ЄЁхї яхЁхьхээ√ї}.
╬э шёяюы№чютры ёыхфє■∙хх яюэ Єшх.
(┼ёыш ¤Єю яюэ Єшх шыш яюёыхфє■∙шщ ЄхъёЄ фю ыхьь√ ╩юыьюуюЁютр юс єэштхЁёры№э√ї
ЇєэъЎш ї яюърцєЄё  ┬рь ёыюцэ√ьш, ┬√ ьюцхЄх ёЁрчє яхЁхщЄш ъ ¤Єющ ыхььх.
─Ёєующ трЁшрэЄ --- яЁюўшЄрЄ№ ¤ЄюЄ ЄхъёЄ фы  $n=2$, р т ¤Єющ ыхььх ёэютр
тхЁэєЄ№ё  ъ яЁюшчтюы№эюьє $n$.)

╬сючэрўшь $I=[-1;1]$.
{\it ─хЁхтюь $T_f$ ъюьяюэхэЄ ьэюцхёЄт єЁютэ } ЇєэъЎшш $f:I^n\to I$
эрч√трхЄё  ьхЄЁшўхёъюх яЁюёЄЁрэёЄтю, Єюўърьш ъюЄюЁюую  ты ■Єё  ъюьяюэхэЄ√
ёт чэюёЄш ьэюцхёЄт $f^{-1}(c)$, $c\in I$, р ьхЄЁшър юяЁхфхыхэр т 
http://en.wikipedia.org/wiki/Hausdorf\_distance

═ряЁшьхЁ, фхЁхтю ъюьяюэхэЄ ьэюцхёЄт єЁютэ  ЇєэъЎшш $f:I^2\to I$,
$f(x,y)=xy$, уюьхюьюЁЇэю сєътх $X$.
─Ёєушх яЁшьхЁ√ яЁштхфхэ√ эр Ёшёєэъх (ёь. фхЄрыш т [Ar58]).

\begin{figure}
  \begin{center}
    \begin{tabular}{@{}c@{\qquad\hspace*{-6pt}}c@{\qquad}c@{\qquad}c@{}}
      \mpfile{plot3d}{1}
      &
      \raisebox{-3.23mm}{\mpfile{plot2d}{1}}
      &
      \mpfile{tree}{2}
      &
      \mpfile{tree}{1}
      \\[3mm]
      уЁрЇшъ ЇєэъЎшш $f$&
      $I^2$\rlap{\hspace*{2mm}%
	$\xrightarrow{\hspace*{10mm}(u_f,v_f)=\bar f\hspace*{10mm}}$}&
      $T_f$\rlap{\hspace*{2mm}$\xrightarrow{\hspace*{3mm}t_f\hspace*{3mm}}$} &
      $I$
      \\[3mm]
    \multicolumn{4}{c}{(a)}\\[6mm]
      \mpfile{plot3d}{2}
      &
    \mpfile{plot2d}{2}
      &
      \mpfile{tree}{3}
      &
      \mpfile{tree}{1}
      \\[3mm]
      уЁрЇшъ ЇєэъЎшш $f$&
      $I^2$\rlap{\hspace*{2mm}%
	$\xrightarrow{\hspace*{10mm}(u_f,v_f)=\bar f\hspace*{10mm}}$}&
      $T_f$\rlap{\hspace*{2mm}$\xrightarrow{\hspace*{3mm}t_f\hspace*{3mm}}$} &
      $I$
      \\[3mm]
    \multicolumn{4}{c}{(b)}
    \end{tabular}
  \end{center}
\caption{╧ЁшьхЁ√ фхЁхт№хт ъюьяюэхэЄ ьэюцхёЄт єЁютэ 
ш Ёрчыюцхэшщ ╩ЁюэЁюфр}\label{pic1}
\end{figure}

╬ўхтшфэю, ўЄю ы■ср  ЇєэъЎш  $f:I^n\to I$ юЄ $n$ яхЁхьхээ√ї
яЁхфёЄрты хЄё  т тшфх ъюьяючшЎшш
$I^n\overset{t_f}\to T_f\overset{\overline f}\to I$ фы  эхъюЄюЁ√ї
юЄюсЁрцхэшщ $\overline f$ ш $t_f$.
╧ЁюёЄЁрэёЄтю $T_f$ ьюцэю ёўшЄрЄ№ ыхцр∙шь схч ёрьюяхЁхёхўхэшщ т ътрфЁрЄх $I^2$.
\footnote{─ы  фюърчрЄхы№ёЄтр эєцэю яюърчрЄ№, ўЄю $T_f$  ты хЄё 
{\it фхЁхтюь}, Є.  х.  юфэюьхЁэ√ь ёЄ уштрхь√ь ыюъры№эю ёт чэ√ь ъюьяръЄюь.
╧ЁшьхЁ√ фхЁхт№хт эрїюф Єё  эр Ёшё. 3,4,5 эшцх. ╦■сюх фхЁхтю яырэрЁэю [Ku68].}
╧ю¤Єюьє $t_f$ ьюцэю ёўшЄрЄ№ ярЁющ ЇєэъЎшщ $u_f,v_f:I^n\to I$.
╘єэъЎш■ $\overline f$ ьюцэю яЁюфюыцшЄ№ эр тхё№ ътрфЁрЄ $I^2$ (яю ЄхюЁхьх
╙Ё√ёюэр ю яЁюфюыцхэшш), Є.х. ёўшЄрЄ№ ЇєэъЎшхщ фтєї яхЁхьхээ√ї.
╚Єръ, шьххь (Ёшё. 1) {\it Ёрчыюцхэшх ╩ЁюэЁюфр}
$$f(x_1,\dots,x_n)=\overline f(u_f(x_1,\dots,x_n),v_f(x_1,\dots,x_n)).$$

{\bf ╦хььр ╩юыьюуюЁютр юс єэштхЁёры№э√ї фхЁхт№ ї.} 
{\it ─ы  ы■сюую $n\ge2$ 

ёє∙хёЄтє■Є Єръшх фхЁхт№  $T_1,\dots,T_{n+1}$ ш ЇєэъЎшш $t_i:I^n\to T_i$, ўЄю

фы  ы■сющ эхяЁхЁ√тэющ ЇєэъЎшш $f:I^n\to I$ юЄ $n$ яхЁхьхээ√ї

ёє∙хёЄтє■Є эхяЁхЁ√тэ√х ЇєэъЎшш $g_1,\dots,g_{n+1}:I^n\to I$ юЄ $n$ яхЁхьхээ√ї,
фы  ъюЄюЁ√ї}
$$T_{g_i}=T_i,\quad t_{g_i}=t_i\quad\mbox{ш}\quad f=g_1+\dots+g_{n+1}.$$

┬рцэю, ўЄю фхЁхт№  $T_{g_i}$ ъюьяюэхэЄ ьэюцхёЄт єЁютэ  ЇєэъЎшщ $g_i$
(ш ёююЄтхЄёЄтє■∙шх юЄюсЁрцхэш  $t_{g_i}$) эх чртшё Є юЄ $f$, їюЄ  ёрьш ЇєэъЎшш
$g_i$ ьюуєЄ чртшёхЄ№ юЄ $f$.

╧ю-тшфшьюьє, ыхььр тхЁэр ш фы  $n=1$ (эю ¤Єю эхЄЁштшры№эю).

─юърчрЄхы№ёЄтр ь√ эх яЁштюфшь.
╒юЄ  юэю  ты ■Єё  трцэ√ь °руюь т фюърчрЄхы№ёЄтх ЄхюЁхь√ ╩юыьюуюЁютр-└Ёэюы№фр,
эю эр°р Ўхы№ --- юётхЄшЄ№ шьхээю Єх °руш, т ъюЄюЁ√ї яю тшыюё№ яюэ Єшх
срчшёэюую тыюцхэш .
╩Ёюьх Єюую, яЁюсыхьє ├шы№схЁЄр ьюцэю Ёх°шЄ№ эрьэюую яЁю∙х: ёь. эшцх
ёєяхЁяючшЎшюээє■ ЄхюЁхьє ╩юыьюуюЁютр ш хх фюърчрЄхы№ёЄтю т [Ar58].

╚ч ыхьь√ ╩юыьюуюЁютр юс єэштхЁёры№э√ї фхЁхт№ ї ш Ёрчыюцхэш  ╩ЁюэЁюфр т√ЄхърхЄ
ёыхфє■∙шщ Ёхчєы№ЄрЄ (фюърцшЄх!).

\smallskip
{\bf ╦хььр ╩юыьюуюЁютр юс єэштхЁёры№э√ї  ЇєэъЎш ї.}
{\it ─ы  ы■сюую $n\ge3$ 

ёє∙хёЄтєхЄ Єръющ эрсюЁ шч $2n+2$ эхяЁхЁ√тэ√ї ЇєэъЎшш $u_i,v_i:I^n\to I$
($i=1,\dots,n+1$) юЄ $n$ яхЁхьхээ√ї,

ўЄю фы  ы■сющ эхяЁхЁ√тэющ ЇєэъЎшш $f:I^n\to I$ юЄ $n$ яхЁхьхээ√ї

ёє∙хёЄтє■Є эхяЁхЁ√тэ√х ЇєэъЎшш $f_i:I^2\to I$ ($i=1,\dots,n+1$) фтєї
яхЁхьхээ√ї, фы  ъюЄюЁ√ї}
$$f(x_1,\dots,x_n)=
\sum\limits_{i=1}^{n+1}f_i(u_i(x_1,\dots,x_n),v_i(x_1,\dots,x_n)).$$

┬рцэю, ўЄю ЇєэъЎшш $u_i,v_i$ эх чртшё Є юЄ $f$ (яЁш ЇшъёшЁютрээюь $n$),
їюЄ  ЇєэъЎшш $f_i$ ьюуєЄ чртшёхЄ№ юЄ $f$.

▌Єр ыхььр ЄЁштшры№эр фы  $n=1$ ш $n=2$ (яюфєьрщЄх, яюўхьє).

\smallskip
{\it ═рсЁюёюъ фюърчрЄхы№ёЄтр ЄхюЁхь√ ╩юыьюуюЁютр ю т√ЁрчшьюёЄш ўхЁхч ЇєэъЎшш
ЄЁхї яхЁхьхээ√ї.}
─ы  ЇєэъЎшш $f(x_1,x_2,x_3,x_4)$ ўхЄ√Ёхї яхЁхьхээ√ї шьххь
$$f(x_1,x_2,x_3,x_4)\ =\ f_{x_4}(x_1,x_2,x_3)\ \overset{(*)}=
\ \sum\limits_{i=1}^4 f_{x_4,i}(u_i(x_1,x_2,x_3),v_i(x_1,x_2,x_3))\ =$$
$$=\ \sum\limits_{i=1}^4 F_i(u_i(x_1,x_2,x_3),v_i(x_1,x_2,x_3),x_4),
\quad\mbox{уфх} \quad F_i(a,b,c)=f_{c,i}(a,b).$$
╘єэъЎш  $f_{x_4}$ эхяЁхЁ√тэю чртшёшЄ юЄ ярЁрьхЄЁр $x_4$.
╨ртхэёЄтю (*) яюыєўрхЄё  яЁшьхэхэшхь ыхьь√ ╩юыьюуюЁютр юс єэштхЁёры№э√ї
ЇєэъЎш ї.
╠√ шёяюы№чєхь єёшыхээ√щ трЁшрэЄ ¤Єющ ыхьь√, єЄтхЁцфр■∙шщ, ўЄю ърцфр  шч $n+1$
ЇєэъЎшщ $f_i$ эхяЁхЁ√тэю чртшёшЄ юЄ шёїюфэющ ЇєэъЎшш $f$.
╚ч ¤Єюую трЁшрэЄр т√ЄхърхЄ эхяЁхЁ√тэр  чртшёшьюёЄ№ ЇєэъЎшщ $f_{x_4,i}$ юЄ
ярЁрьхЄЁр $x_4$ ($i=1,2,3,4$).
└ ¤Єю тыхўхЄ эхяЁхЁ√тэюёЄ№ ЇєэъЎшщ $F_i$ ($i=1,2,3,4$).

─ы  ЇєэъЎшш сюы№°хую ъюышўхёЄтр яхЁхьхээ√ї Ёрёёєцфхэшх рэрыюушўэю. QED

\smallskip
{\bf 1.} ╟р юфэє ъюяхщъє ртЄюьрЄ т√фрхЄ чэрўхэшх чрфрээющ ┬рьш эхяЁхЁ√тэющ
ЇєэъЎшш ЄЁхї яхЁхьхээ√ї эр чрфрээющ ┬рьш ЄЁющъх ўшёхы.
╟р ъръє■ ёєььє ┬√ чртхфюью ёьюцхЄх т√ўшёышЄ№ чрфрээє■ эхяЁхЁ√тэє■ ЇєэъЎш■ $n$
яхЁхьхээ√ї (яЁш єёыютшш эрышўш  є ┬рё эхюуЁрэшўхээющ ярь Єш)?

\bigskip
{\bf ╥хюЁхьр └Ёэюы№фр: ъ ёєяхЁяючшЎш ь ЇєэъЎшщ фтєї яхЁхьхээ√ї.}

─ы  фюърчрЄхы№ёЄтр ЄхюЁхь√ ╩юыьюуюЁютр-└Ёэюы№фр юёЄрыюё№ яЁюшчтюы№эє■
эхяЁхЁ√тэє■ ЇєэъЎш■ ЄЁхї яхЁхьхээ√ї т√ЁрчшЄ№ ўхЁхч
ёєяхЁяючшЎш■ эхяЁхЁ√тэ√ї ЇєэъЎшщ фтєї яхЁхьхээ√ї.
─ы  ¤Єюую яюыхчэю ёыхфє■∙хх яюэ Єшх.

╧юфьэюцхёЄтю $T\subset I^3$ эрчютхь {\it срчшёэ√ь}, хёыш ы■ср  эхяЁхЁ√тэр 
ЇєэъЎш  эр $T$ ьюцхЄ с√Є№ Ёрчыюцхэр т ёєььє ЄЁхї ЇєэъЎшщ, ърцфр  шч
ъюЄюЁ√ї чртшёшЄ Єюы№ъю юЄ юфэющ ъююЁфшэрЄ√. ╚ыш, ЇюЁьры№эю,
хёыш фы  ы■сющ эхяЁхЁ√тэющ ЇєэъЎшш $f:T\to I$ ёє∙хёЄтє■Є Єръшх эхяЁхЁ√тэ√х
ЇєэъЎшш $g_1,g_2,g_3:I\to I$, ўЄю
$$f(x,y,z)=g_1(x)+g_2(y)+g_3(z)\quad\mbox{фы }\quad (x,y,z)\in T.$$

{\bf ╦хььр └Ёэюы№фр ю фхЁхт№ ї.}
{\it ╦■сюх фхЁхтю ЁхрышчєхЄё  ъръ срчшёэюх яюфьэюцхёЄтю т $I^3$
(Є.х. Єюяюыюушўхёъш ¤ътштрыхэЄэю эхъюЄюЁюьє срчшёэюьє яюфьэюцхёЄтє т $I^3$).
\footnote{─юърчрЄхы№ёЄтю ыхьь√ └Ёэюы№фр шёяюы№чєхЄ ЄхюЁхьє ╠хэухЁр ю
ёє∙хёЄтютрэшш єэштхЁёры№эюую фхЁхтр.
═р ёрьюь фхых, └Ёэюы№ф фюърчры ¤Єє ыхььє фы  фхЁхт№хт ё Єюўърьш
тхЄтыхэш  ЄЁхЄ№хую яюЁ фър.
▌Єюую с√ыю фюёЄрЄюўэю фы  Ёх°хэш  яЁюсыхь√ ├шы№схЁЄр.
╬с∙шщ ёыєўрщ ыхьь√ фюърчрэ ╬ёЄЁрэфюь т 1965 [St89].}
}

\smallskip
╚фх  фюърчрЄхы№ёЄтр тшфэр эр яЁшьхЁх фюърчрЄхы№ёЄтр ышсю срчшёэющ тыюцшьюёЄш т
яыюёъюёЄ№ ъюэхўэюую фхЁхтр, шч ърцфющ тхЁ°шэ√ ъюЄюЁюую т√їюфшЄ ышсю юфэю, ышсю
ЄЁш ЁхсЁр [Ar58], ышсю сюыхх ёшы№эюую єЄтхЁцфхэш  (c) т яЁхфяюёыхфэхь яєэъЄх 
тЄюЁющ ўрёЄш.

\smallskip
{\bf ╦хььр └Ёэюы№фр юс єэштхЁёры№э√ї ЇєэъЎш ї.}
{\it ╤є∙хёЄтєхЄ Єръющ эрсюЁ шч фхт Єш эхяЁхЁ√тэ√ї ЇєэъЎшш
$u_i:I^2\to I$ ($i=1,2,\dots,9$) фтєї яхЁхьхээ√ї, ўЄю

фы  ы■сющ эхяЁхЁ√тэющ ЇєэъЎшш $f:I^2\to I$ фтєї яхЁхьхээ√ї

ёє∙хёЄтє■Є эхяЁхЁ√тэ√х ЇєэъЎшш
$f_i:I\to I$ ($i=1,2,\dots,9$) юфэющ яхЁхьхээющ, 

фы  ъюЄюЁ√ї
$f(x,y)=f_1(u_1(x,y))+\dots+f_9(u_9(x,y))$.}

\smallskip
┬рцэю, ўЄю ЇєэъЎшш $u_i$ эх чртшё Є юЄ $f$, їюЄ  ЇєэъЎшш $f_i$ ьюуєЄ чртшёхЄ№ 
юЄ $f$.

\smallskip
{\it ─юърчрЄхы№ёЄтю.}
┬юч№ьхь фхЁхт№  $T_1,T_2,T_3$ шч ыхьь√ ╩юыьюуюЁютр юс єэштхЁёры№э√ї фхЁхт№ ї 
фы  $n=2$.
╧ю ыхььх └Ёэюы№фр ю фхЁхт№ ї ёє∙хёЄтє■Є срчшёэ√х тыюцхэш 
$(u_{i1},u_{i2},u_{i3}):T_i\to I^3$.
╧юыюцшь $u_{3(i-1)+j}:=u_{ij}$. 
┬юч№ьхь ЇєэъЎшш $g_1,g_2,g_3:I^2\to I$ (чртшё ∙шх юЄ $f$) шч ыхьь√ ╩юыьюуюЁютр 
юс єэштхЁёры№э√ї фхЁхт№ ї фы  $n=2$.
╧Ёшьхэ   рэрыюу Ёрчыюцхэш  ╩ЁюэЁюфр ш юяЁхфхыхэшх срчшёэюёЄш, яюыєўрхь  
$$g_i(x,y)=\overline g_i(u_{i1}(x,y),u_{i2}(x,y),u_{i3}(x,y))=
g_{i1}(u_{i1}(x,y))+g_{i2}(u_{i2}(x,y))+g_{i3}(u_{i3}(x,y))$$
фы  эхъюЄюЁ√ї ЇєэъЎшщ $g_{ij}:I^2\to I$. 
╬ёЄрхЄё  яюыюцшЄ№ $f_{3(i-1)+j}:=g_{ij}$. QED

\smallskip
╥хяхЁ№ ЄхюЁхьр ╩юыьюуюЁютр-└Ёэюы№фр т√ЄхърхЄ шч
$$f(x,y,z)=f_z(x,y)=\sum\limits_{i=1}^9 f_{i,z}(u_i(x,y))=
\sum\limits_{i=1}^9 F_i(u_i(x,y),z),\quad\mbox{уфх}
\quad F_i(t,z)=f_{i,z}(t).$$
═хяЁхЁ√тэюёЄ№ ЇєэъЎшщ $F_i$ фюърч√трхЄё  рэрыюушўэю яЁхф√фє∙хьє яєэъЄє  
(ё шёяюы№чютрэшхь ёююЄтхЄёЄтє■∙хую єёшыхэш  ыхьь√ └Ёэюы№фр юс єэштхЁёры№э√ї 
ЇєэъЎш ї).

\smallskip
{\bf 2.} ╟р юфэє ъюяхщъє ртЄюьрЄ т√фрхЄ чэрўхэшх чрфрээющ ┬рьш эхяЁхЁ√тэющ
ЇєэъЎшш
фтєї яхЁхьхээ√ї эр чрфрээющ ┬рьш ярЁх ўшёхы.
╟р ъръє■ ёєььє ┬√ чртхфюью ёьюцхЄх т√ўшёышЄ№ чрфрээє■ эхяЁхЁ√тэє■ ЇєэъЎш■ $n$
яхЁхьхээ√ї (яЁш эрышўшш эхюуЁрэшўхээющ ярь Єш)?

\bigskip
{\bf ╥хюЁхьр ╩юыьюуюЁютр: ъ ЇєэъЎш ь юфэющ яхЁхьхээющ ш ёыюцхэш■.}

┬ Єюь цх 1957 уюфє ╩юыьюуюЁют фюърчры х∙х сюыхх ёшы№э√щ Ёхчєы№ЄрЄ, шч ъюЄюЁюую
Єръцх т√ЄхърхЄ Ёх°хэшх яЁюсыхь√ ├шы№схЁЄр.

\smallskip
{\bf ╤єяхЁяючшЎшюээр  ЄхюЁхьр ╩юыьюуюЁютр.}
{\it ╦■ср  эхяЁхЁ√тэр  ЇєэъЎш  яЁхфёЄрты хЄё  т тшфх ёєяхЁяючшЎшш ёыюцхэш  ш
эхяЁхЁ√тэ√ї ЇєэъЎшщ юфэющ яхЁхьхээющ.

┴юыхх Єюўэю, фы  ърцфюую $n>1$ ёє∙хёЄтєхЄ эрсюЁ $n(2n+1)$ Єръшї эхяЁхЁ√тэ√ї
ЇєэъЎшщ $u_{ij}:I\to I$ ($i=1,\dots,2n+1$, $j=1,\dots,n$) юфэющ яхЁхьхээющ,
ўЄю

фы  ы■сющ эхяЁхЁ√тэющ ЇєэъЎшш $f:I^n\to I$ юЄ $n$ яхЁхьхээ√ї

ёє∙хёЄтє■Є
эхяЁхЁ√тэ√х ЇєэъЎшш $f_1,\dots,f_{2n+1}:I\to I$ юфэющ яхЁхьхээющ, фы  ъюЄюЁ√ї}
$$f(x_1,\dots,x_n)\ =\ \sum\limits_{i=1}^{2n+1}
f_i\left(\sum\limits_{j=1}^nu_{ij}(x_j)\right).$$

╟фхё№ трцэю, ўЄю ЇєэъЎшш $u_{ij}$ эх чртшё Є юЄ $f$, їюЄ  ЇєэъЎшш $f_i$
ьюуєЄ чртшёхЄ№ юЄ $f$.
▌ыхьхэЄрЁэюх шчыюцхэшх фюърчрЄхы№ёЄтр яЁштхфхээющ ЄхюЁхь√ ьюцэю эрщЄш т [Ar58].

\smallskip
{\bf 3.} ╟р юфэє ъюяхщъє ртЄюьрЄ ышсю ёъырф√трхЄ фтр чрфрээ√ї ┬рьш ўшёыр,
ышсю т√фрхЄ чэрўхэшх чрфрээющ ┬рьш эхяЁхЁ√тэющ ЇєэъЎшш юфэющ яхЁхьхээющ т
чрфрээющ ┬рьш Єюўъх.
╟р ъръє■ ёєььє ┬√ чртхфюью ёьюцхЄх т√ўшёышЄ№ чрфрээє■ эхяЁхЁ√тэє■ ЇєэъЎш■ $n$
яхЁхьхээ√ї (яЁш эрышўшш эхюуЁрэшўхээющ ярь Єш)?


\smallskip
╬с рэрышЄшўхёъшї яЁюсыхьрї, ёт чрээ√ї ё ¤Єшь т√фр■∙шьё  Ёхчєы№ЄрЄюь
╩юыьюуюЁютр, ёь. [St89, Vi04].
╬ Єюяюыюушўхёъшї яЁюсыхьрї эряшёрэю фрыхх.

\bigskip
{\bf ┴рчшёэ√х тыюцхэш  т ьэюуюьхЁэ√х яЁюёЄЁрэёЄтр.}
\footnote{▌ЄюЄ яєэъЄ эх¤ыхьхэЄрЁхэ, ЇюЁьры№эю эх шёяюы№чєхЄё  т фры№эхщ°хь ш
ьюцхЄ с√Є№ юяє∙хэ ўшЄрЄхыхь.
╬фэръю ь√ яЁштюфшь хую, яюёъюы№ъє юэ фрхЄ ўхЄъє■ ърЁЄшэє фы  Ёрчэ√ї
ЁрчьхЁэюёЄхщ.}

╧юфьэюцхёЄтю $K\subset\R^n$ эрч√трхЄё  {\it срчшёэ√ь}, хёыш фы  ы■сющ
эхяЁхЁ√тэющ ЇєэъЎшш $f:K\to\R$ ёє∙хёЄтє■Є Єръшх эхяЁхЁ√тэ√х ЇєэъЎшш
$f_1,\dots,f_n:\R\to\R$, ўЄю
$$f(x_1,\dots,x_n)=f_1(x_1)+\dots+f_n(x_n)\quad\mbox{фы }\quad
(x_1,\dots,x_n)\in K.$$
╧ЁюёЄЁрэёЄтю $K$ эрч√трхЄё  {\it срчшёэю тыюцшь√ь} т $\R^n$, хёыш
ёє∙хёЄтєхЄ тыюцхэшх $K\to\R^n$, юсЁрч ъюЄюЁюую срчшёэ√щ.

╘єэъЎшш эр яЁюшчтюы№эюь $n$-ьхЁэюь ъюьяръЄх єцх эхы№ч  яЁхфёЄрты Є№ ёхсх
ъръ ЇєэъЎшш $n$ яхЁхьхээ√ї.
╬фэръю яюэ Єшх срчшёэющ тыюцшьюёЄш фюёЄрты хЄ рэрыюу ЁрчыюцшьюёЄш ЇєэъЎшщ
эр ъюьяръЄрї т ёєяхЁяючшЎш■ ЇшъёшЁютрээ√ї ЇєэъЎшщ ш ёыюцхэш .

╚ч ёєяхЁяючшЎшюээющ ЄхюЁхь√ ╩юыьюуюЁютр ёыхфєхЄ, ўЄю {\it $n$-ьхЁэ√щ ъєс
срчшёэю тыюцшь т $\R^{2n+1}$.}
─хщёЄтшЄхы№эю, ЇєэъЎшш $u_i=\sum\limits_{j=1}^n u_{ij}$ ($i=1,\dots,2n+1$)
шч ЄхюЁхь√ ╩юыьюуюЁютр юяЁхфхы ■Є срчшёэюх тыюцхэшх $I^n\to I^{2n+1}$.
╬ёЄЁрэф чрьхЄшы т 1965~у., ўЄю ¤ЄюЄ ЇръЄ ьюцэю юсюс∙шЄ№.

\smallskip
{\bf ╥хюЁхьр ╬ёЄЁрэфр.}
{\it ╦■сющ $n$-ьхЁэ√щ ъюьяръЄ срчшёэю тыюцшь т $\R^{2n+1}$} [St89].

\smallskip
▌Єр ЄхюЁхьр єёшыштрхЄ ЄхюЁхьє ═хсышэур--╠хэухЁр--╧юэЄЁ ушэр ю
тыюцшьюёЄш ы■сюую $n$-ьхЁэюую ъюьяръЄр т $\R^{2n+1}$ [Ku68].

═р ёрьюь фхых ╬ёЄЁрэф фюърчры ёыхфє■∙шщ сюыхх ёшы№э√щ Ёхчєы№ЄрЄ, юсюс∙р■∙шщ
ёєяхЁяючшЎшюээє■ ЄхюЁхьє ╩юыьюуюЁютр (р эх Єюы№ъю хх ёыхфёЄтшх).

{\it  ╧єёЄ№ $X_1,\dots,X_m$ --- ъюэхўэюьхЁэ√х ьхЄЁшўхёъшх яЁюёЄЁрэёЄтр.
╧юыюцшь $n=\dim X_1+\dots+\dim X_m$ ш $X=X_1\times\dots\times X_m$.
╥юуфр ёє∙хёЄтє■Є Єръшх эхяЁхЁ√тэ√х ЇєэъЎшш $u_{ij}:X_j\to\R$,
($i=1,\dots,2n+1$, $j=1,\dots,m$), ўЄю фы  ЇєэъЎшщ
$u_i(x_1,\dots,x_m)=u_{i1}(x_1)+\dots+u_{im}(x_m)$ ш ы■сющ эхяЁхЁ√тэющ ЇєэъЎшш
$f:X\to\R$ ёє∙хёЄтє■Є эхяЁхЁ√тэ√х ЇєэъЎшш $f_1,\dots,f_{2n+1}:\R\to\R$, фы 
ъюЄюЁ√ї}
$$f(x_1,\dots,x_m)=f_1(u_1(x_1,\dots x_m))+\dots+
f_{2n+1}(u_{2n+1}(x_1,\dots,x_m)).$$

╚ьх■Єё  $n$-ьхЁэ√х яюыш¤фЁ√, эх тыюцшь√х т $\R^{2n}$ [Pr04, Sk].

\smallskip
{\bf ╥хюЁхьр ╪ЄхЁэЇхы№фр.} {\it ─ы  ы■сюую $n\ge2$ ы■сющ $n$-ьхЁэ√щ
ъюьяръЄ (эряЁшьхЁ, $n$-ьхЁэ√щ ъєс) эх тыюцшь срчшёэю т $\R^{2n}$} [St89].


\smallskip
╚эЄхЁхёэю, ўЄю ЄхюЁхьр ╪ЄхЁэЇхы№фр ЁхфєЎшЁєхЄё  ъ ъюьсшэрЄюЁэю-ухюьхЄЁшўхёъюьє
єЄтхЁцфхэш■ яЁш яюью∙ш ьэюуюьхЁэюую рэрыюур ъЁшЄхЁш  срчшёэюёЄш, яЁштхфхээюую
т ЄЁхЄ№хщ ўрёЄш эрёЄю ∙хщ ёЄрЄ№ш.

╬ўхтшфэю, ўЄю $K$ срчшёэю тыюцшь т $\R$ Єюуфр ш Єюы№ъю Єюуфр, ъюуфр $K$
Єюяюыюушўхёъш тыюцшь т $\R$.
╚ч ЄхюЁхь ╬ёЄЁрэфр ш ╪ЄхЁэЇхы№фр ёыхфєхЄ, ўЄю {\it фы  $m>2$
ъюьяръЄ $K$ срчшёэю тыюцшь т $\R^m$ Єюуфр ш Єюы№ъю Єюуфр, ъюуфр
$\dim K<m/2$.}
╥ръшь юсЁрчюь, т 1989~у. юёЄртрыюё№ эхшчтхёЄэ√ь ыш°№ юяшёрэшх ъюьяръЄют, срчшёэю
тыюцшь√ї {\it т яыюёъюёЄ№}.

\bigskip
\newpage
\centerline{\uppercase{\bf ┴рчшёэ√х тыюцхэш  т яыюёъюёЄ№}}

\smallskip
{\bf ┴рчшёэр  тыюцшьюёЄ№ т яыюёъюёЄ№.}

├ЁрЇ (шыш ъюьяръЄ) $K$ эрч√трхЄё  {\it срчшёэю тыюцшь√ь} т яыюёъюёЄ№, хёыш
ёє∙хёЄтєхЄ Єръюх тыюцхэшх $\varphi:K\to\R^2$, ўЄю фы  ы■сющ эхяЁхЁ√тэющ ЇєэъЎшш
$f:\varphi(K)\to\R$ ёє∙хёЄтє■Є Єръшх эхяЁхЁ√тэ√х ЇєэъЎшш
$g,h:\R\to\R$, ўЄю $f(x,y)=g(x)+h(y)$ фы  ы■сющ Єюўъш $(x,y)\in\varphi(K)$.
(╬яЁхфхыхэшх эхяЁхЁ√тэющ ЇєэъЎшш эряюьэхэю т эрўрых ўрёЄш 1.)

╧Ёюсыхьр юяшёрэш  уЁрЇют (ш ъюьяръЄют), срчшёэю тыюцшь√ї т яыюёъюёЄ№,
яюёЄртыхэр ╪ЄхЁэЇхы№фюь [St89].
╩ЁшЄхЁшщ срчшёэющ тыюцшьюёЄш ышэхщэю-ёт чэ√ї ъюьяръЄют т яыюёъюёЄ№ яюыєўхэ т
[Sk95].
─ы  ъюэхўэ√ї уЁрЇют юэ ЇюЁьєышЁєхЄё  юёюсхээю яЁюёЄю.

\begin{figure}
  \begin{center}
    \begin{tabular}{c@{\qquad}c@{\qquad}c}
      \mpfile{gr}{1}
      &
      \mpfile{gr}{2}
      &
      \mpfile{gr}{3}
      \\[3mm]
      $S^1$& $T_5$ & $C$
    \end{tabular}
  \end{center}
  \caption{}\label{pic2}
\end{figure}

\begin{figure}
  \begin{center}
    \begin{tabular}{c@{\qquad}c@{\qquad}c}
      \raisebox{12.9mm}{\mpfile{gr}{31}}
      &
      \raisebox{3.87mm}{\mpfile{gr}{32}}
      &
      \mpfile{gr}{33}
      \\[3mm]
      $R_1$& $R_2$ & $R_3$\\[5mm]
      \raisebox{12.9mm}{\mpfile{gr}{34}}
      &
      \raisebox{3.87mm}{\mpfile{gr}{35}}
      &
      \mpfile{gr}{36}
      \\[3mm]
      $F_1$& $F_2$ & $F_3$
    \end{tabular}
  \end{center}
  \caption{}\label{pic3}
\end{figure}


\smallskip
{\bf ╩ЁшЄхЁшщ срчшёэющ тыюцшьюёЄш уЁрЇют.} [Sk95, ёЁ. Sk05]
{\it ╩юэхўэ√щ уЁрЇ $K$ срчшёэю тыюцшь т яыюёъюёЄ№ Єюуфр ш Єюы№ъю Єюуфр, ъюуфр
т√яюыэхэю юфэю шч ёыхфє■∙шї фтєї ¤ътштрыхэЄэ√ї єёыютшщ:

(S) $K$ эх ёюфхЁцшЄ яюфуЁрЇют, уюьхюьюЁЇэ√ї юъЁєцэюёЄш $S^1$, яхэЄюфє
$T_5=C_1$ шыш ъЁхёЄє ё ЁрчтхЄтыхээ√ьш ъюэЎрьш $C=C_2$ (Ёшё. 2);

(U) $K$ ёюфхЁцшЄё  т юфэюь шч уЁрЇют $R_n$ (Ёшё. 3).}

\begin{figure}
  \begin{center}
    \begin{tabular}{c@{\qquad}c@{\qquad}c}
      \mpfile{gr}{41}
      &
      \mpfile{gr}{42}
      &
      \mpfile{gr}{43}
      \\[3mm]
      $C_4$& $B$ & $H_2$
    \end{tabular}
  \end{center}
  \caption{}\label{pic4}
\end{figure}

\begin{figure}
  \begin{center}
    \begin{tabular}{c@{\qquad}c}
      \mpfile{gr}{51}
      &
      \mpfile{gr}{52}
      \\[2mm]
      $C_3$& $F$
      \\[3mm]
      \mpfile{gr}{53}
      & 
      \mpfile{gr}{54}
      \\[2mm]
      $H_-$&$H_+$
      \\[3mm]
      \mpfile{gr}{55}
       &
      \mpfile{gr}{56}
      \\[2mm]
      $h_+$& $h_-$  
    \end{tabular}
  \end{center}
  \caption{}\label{pic5}
\end{figure}

\smallskip
{\it ╬яЁхфхыхэшх уЁрЇют $F_n$ ш $R_n$ (Ёшё. 3).}
╧єёЄ№ $F_1$~--- ЄЁшюф, ш фы  ы■сюую $n$ уЁрЇ $F_{n+1}$ яюыєўхэ шч $F_n$
ЁрчтхЄтыхэшхь ърцфюую тшё ўхую ЁхсЁр уЁрЇр $F_n$.
├ЁрЇ $R_n$ яюыєўрхЄё  шч уЁрЇр $F_n$ фюсртыхэшхь тшё ўхую ЁхсЁр ъ ърцфющ Єюўъх
тхЄтыхэш  уЁрЇр $F_n$.

\smallskip
─юърчрЄхы№ёЄтю яЁштюфшЄё  т ёыхфє■∙хь яєэъЄх.

╧Ёштхфхь схч фюърчрЄхы№ёЄтр Ёх°хэшх яЁюсыхь√ ╪ЄхЁэЇхы№фр фы  сюыхх юс∙хую
ёыєўр 
(хую ьюцэю яЁюяєёЄшЄ№ схч є∙хЁср фы  яюэшьрэш  фры№эхщ°хую).

\smallskip
{\bf ╩ЁшЄхЁшщ срчшёэющ тыюцшьюёЄш ышэхщэю-ёт чэ√ї ъюьяръЄют.} [Sk95, ёЁ. Sk05]
{\it ╦шэхщэю-ёт чэ√щ ъюьяръЄ $K$ срчшёэю тыюцшь т яыюёъюёЄ№ Єюуфр ш Єюы№ъю
Єюуфр, ъюуфр юэ  ты хЄё  ыюъры№эю ёт чэ√ь (Є.х. яхрэютёъшь) ш т√яюыэхэю юфэю шч
ёыхфє■∙шї фтєї ¤ътштрыхэЄэ√ї єёыютшщ:

(1) $K$ эх ёюфхЁцшЄ яюфъюьяръЄют $S^1, C_2, C_4, B$ ш, фы  эхъюЄюЁюую $n$,
яюфъюьяръЄют $F_n$ ш $H_n$ (Ёшё. 2, 3, 4);

(2) $K$ эх ёюфхЁцшЄ яюфъюэЄшэєєьют
$S^1, C_1, C_2, C_3, B, F, H_+$, $H_-, h_+, h_-$ (Ёшё. 2, 4, 5).}


\smallskip
┬тхфхь шёяюы№чютрээ√х юсючэрўхэш  (Ёшё. 4, 5).
{\it ═єы№-яюёыхфютрЄхы№эюёЄ№■} ьэюцхёЄт эрч√трхЄё  яюёыхфютрЄхы№эюёЄ№
ьэюцхёЄт, фшрьхЄЁ√ ъюЄюЁ√ї ёЄЁхь Єё  ъ эєы■.

╬сючэрўшь ўхЁхч $C_3$~ъЁхёЄ ё эєы№-яюёыхфютрЄхы№эюёЄ№■ фєу, ёїюф ∙шїё  ъ хую
'ЎхэЄЁє' ш яЁшъыххээ√ї ъ юфэющ шч хую 'тхЄтхщ'.

╬сючэрўшь ўхЁхч $C_4$~ъЁхёЄ ё яюёыхфютрЄхы№эюёЄ№■ Єюўхъ, ёїюф ∙шїё  ъ хую
'ЎхэЄЁє'.

╩рцф√щ шч ъюэЄшэєєьют $B,H_n,F,H_+,H_-$  ты хЄё  юс·хфшэхэшхь юЄЁхчър $I=[0;1]$
ш эєы№-яюёыхфютрЄхы№эюёЄш

$\bullet$ фєу, яЁшъыххээ√ї чр юфшэ ъюэхЎ ъ $(0,1)$ т фтюшўэю-ЁрЎшюэры№э√ї
Єюўърї (фы  $B$; юўхтшфэю, ўЄю Єюяюыюушўхёъшщ Єшя яЁюёЄЁрэёЄтр $B$
эх чртшёшЄ юЄ трЁшрЎшщ т ¤Єюь яюёЄЁюхэшш);

$\bullet$ ЄЁшюфют, яЁшъыххээ√ї ъ $I$ чр юфшэ ъюэхЎ т Єюўърї
ьэюцхёЄтр $\{3^{-l_1}+\dots+3^{-l_s}\ |\ s\le n,\
0<l_1<\dots<l_s\mbox{ -- Ўхы√х}\}$ (фы  $H_n$);

$\bullet$ уЁрЇют $F_n$, яЁшъыххээ√ї ъ Єюўърь $1/n$ чр юфэє шч тшё ўшї
тхЁ°шэ (фы  $F$);

$\bullet$ ъюэЄшэєєьют $H_n$, ёюхфшэхээ√ї ё Єюўърьш $1/n$ фєурьш,
яхЁхёхър■∙шьш $H_n$ т $1\in I\subset H_n$ (фы  $H_+$) шыш т $0\in I\subset H_n$
(фы  $H_-$).

╩юэЄшэєєь $h_+$ (ёююЄтхЄёЄтхээю $h_-$) яюыєўхэ шч эєы№-яюёыхфютрЄхы№эюёЄш
ъюэЄшэєєьют $H_n$ ёъыхштрэшхь Єюўхъ $1\in I\subset H_n$ ш
$0\in I\subset H_{n-1}$ (ёююЄтхЄёЄтхээю Єюўхъ $0\in I\subset H_n$ ш
$1\in I\subset H_{n-1}$).


├шяюЄхчр ю срчшёэющ тыюцшьюёЄш (эх юс чрЄхы№эю ышэхщэю-ёт чэ√ї) ъюэЄшэєєьют т
яыюёъюёЄ№ х∙х сюыхх уЁюьючфър.
╬эр ёЇюЁьєышЁютрэр т [Sk95, RS99].

\bigskip
{\bf ═рсЁюёюъ фюърчрЄхы№ёЄтр ъЁшЄхЁш  срчшёэющ тыюцшьюёЄш уЁрЇют.}

─юёЄрЄюўэю фюърчрЄ№ ёыхфє■∙шх ЄЁш єЄтхЁцфхэш .

(a) ╬ъЁєцэюёЄ№ $S^1$, яхэЄюф $T_5$ ш ъЁхёЄ ё ЁрчтхЄтыхээ√ьш ъюэЎрьш $C$
(Ёшё. 2) эх тыюцшь√ срчшёэю т $\R^2$.

(b) ┼ёыш ъюэхўэ√щ уЁрЇ $K$ эх ёюфхЁцшЄ эш юфэюую шч уЁрЇют $S^1$, $T_5$ ш $C$
(Ёшё. 2), Єю $K$ ёюфхЁцшЄё  т $R_n$ (Ёшё. 3) фы  эхъюЄюЁюую $n$.

(c) ╩рцф√щ уЁрЇ $R_n$ (Ёшё. 3) срчшёэю тыюцшь т $\R^2$.

\smallskip
{\it ═рсЁюёюъ фюърчрЄхы№ёЄтр єЄтхЁцфхэш  (b).}
 ─юърцхь, ўЄю {\it фхЁхтю $K$ ё $n$ эхтшё ўшьш тхЁ°шэрьш, эх ёюфхЁцр∙хх уЁрЇют
$T_5$ ш $C$, ёюфхЁцшЄё  т $R_n$}.
 ┬юч№ьхь яЁюшчтюы№эє■ тхЁ°шэє $a\in K$.
 ╧юёъюы№ъє $K$ эх ёюфхЁцшЄ $T_5$ ш $C$, Єю $\deg a\leq4$,
  яЁшўхь хёыш $\deg a=4$, Єю $a$ шьххЄ тшё ўхх ЁхсЁю.
 ╟эрўшЄ, юъЁхёЄэюёЄ№ Єюўъш $a$ шч $K$ ьюцэю тыюцшЄ№ т $R_n$ Єръ, ўЄю $a$
  яюярфхЄ т ЎхэЄЁ $R_n$, р ¤Єр юъЁхёЄэюёЄ№ яхЁхщфхЄ т юъЁхёЄэюёЄ№ $T_4$
  ЎхэЄЁр $R_n$.
 ╤ ърцфющ тхЁ°шэющ, ёюёхфэхщ ё $a$, яюёЄєярхь рэрыюушўэю.
 ╧юёъюы№ъє 'уыєсшэр' уЁрЇр $R_n$ (Є.х. ъюышўхёЄтю эхтшё ўшї тхЁ°шэ эр ёрьющ
 фышээющ тхЄтш юЄ ЎхэЄЁр) Ёртэр $n$, р эхтшё ўшї тхЁ°шэ т $K$ Ёютэю $n$,
  Єю, яЁюфюыцр  ¤ЄюЄ яЁюЎхёё фры№°х, ь√ ёьюцхь тыюцшЄ№ т $R_n$ тхё№ уЁрЇ
  $K$. QED



\begin{figure}
  \begin{center}
    \begin{tabular}{@{}c@{\qquad}c@{\qquad}c@{}}
      \mpfile{embed}{61}
      &       &
      \mpfile{embed}{62}\\[3mm]
      (a)&&(b)
    \end{tabular}
  \end{center}
  \caption{}\label{pic6}
\end{figure}

\smallskip
{\it ═рсЁюёюъ фюърчрЄхы№ёЄтр єЄтхЁцфхэш  (ё).}
╬сючэрўшь ўхЁхч $R_0$ юЄЁхчюъ.
┴рчшёэ√х тыюцхэш  $R_n\to\R^2$ ёЄЁю Єё  яю шэфєъЎшш фы  $n\ge0$.
┬ыюцшь $R_0$ т $[-7;5]^2$ ъръ фшруюэры№, ёюхфшэ ■∙є■ Єюўъш $(-7,-7)$ ш 
$(5,5)$.

┬ыюцхэшх $R_n\to\R^2$ яюыєўрхЄё  шч тыюцхэш  эр Ёшё.6a фюсртыхэшхь ЄЁхї 
тыюцхэшщ уЁрЇр $R_{n-1}$ т ътрфЁрЄшъш $A$, $B$, $C$.
╧ЁюхъЎшш $A_x,B_x,C_x$ ътрфЁрЄшъют эр юё№ $Ox$ эх яхЁхёхър■Єё  фЁєу ё фЁєуюь ш 
ё яЁюхъЎшхщ юЄЁхчър т $R_n$, ярЁрыыхы№эюую уюЁшчюэЄры№эющ юёш.
└эрыюушўэюх єЄтхЁцфхэшх ёяЁртхфыштю ш фы  яЁюхъЎшщ $A_y,B_y,C_y$ ътрфЁрЄшъют 
эр юё№ $Oy$.  

─юърцхь, ўЄю яюёЄЁюхээ√х тыюцхэш  срчшёэ√х, яЁш яюью∙ш шэфєъЎшш яю $n$.
┴рчр шэфєъЎшш $n=0$ юўхтшфэр.
─юърцхь °ру шэфєъЎшш.


╧єёЄ№ $n\ge1$ ш $f:R_n\to\R$ --- эхяЁхЁ√тэр  ЇєэъЎш .
─ы  $t\in[0,1]$ яюыюцшь $g(t):=f(t,0)$,
$$h(t):=f(t,t)-g(t)=f(t,t)-f(t,0)\quad\text{ш}
\quad g(-t):=f(-t,t)-h(t)=f(-t,t)-f(t,t)+f(t,0).$$
╧ю яЁхфяюыюцхэш■ шэфєъЎшш ёє∙хёЄтє■Є Єръшх ЇєэъЎшш
$$g:A_x\cup B_x\cup C_x\to\R\quad\text{ш}\quad h:A_y\cup B_y\cup C_y\to\R,
\quad\text{ўЄю}$$
$$f(x,y)=g(x)+h(y)\quad\text{фы }\quad
(x,y)\in (A\cup B\cup C)\cap R_n.$$
╧юыюцшь
$$g(-t):=f(-t,t)-h(t)\quad\text{фы }\quad t\in C_y
\quad\text{ш}\quad h(t)=f(t,t)-g(t)\quad\text{фы }
\quad t\in(-C_y)\cup B_x\to\R.$$
╧Ёюфюыцшь яюёЄЁюхээє■ ЇєэъЎш■ $g:A_x\cup B_x\cup (-C_y)\cup[-1;1]\cup C_x\to\R$
фю эхяЁхЁ√тэющ ЇєэъЎшш $g:[-7;1]\cup C_x\to\R$ (эряЁшьхЁ, ышэхщэю).
╧юыюцшь
$$h(t)=f(-|t|,t)-g(t)\quad\text{фы }\quad t\in[-6;4]\quad\text{ш}\quad
g(t):=f(t,t)-h(t)\quad\text{фы }\quad t\in[1;2]\cup[3;5]$$
(¤Єю юяЁхфхыхэшх ёютярфрхЄ c яЁхцэшь фы  $t\in[-1;1]$).
╧юёых ¤Єюую яЁюфюыцшь $g$ ш $h$ фю эхяЁхЁ√тэ√ї ЇєэъЎшщ $\R\to\R$.
▀ёэю, ўЄю яюыєўхээ√х ЇєэъЎшш $g$ ш $h$ --- шёъюь√х.
QED

\smallskip
╧Ёш фюърчрЄхы№ёЄтх єЄтхЁцфхэш  (a) ь√ шёяюы№чєхь ъЁшЄхЁшщ срчшёэюёЄш шч
ўрёЄш 3 (р Єръцх яЁштхфхээюх яхЁхф эшь юяЁхфхыхэшх ьюыэшш ш яЁштхфхээюх яюёых
эхую юяЁхфхыхэшх юяхЁрЎшш $E$).

\smallskip
{\it ─юърчрЄхы№ёЄтю срчшёэющ эхтыюцшьюёЄш юъЁєцэюёЄш.} [St89]
╧єёЄ№ чрфрэю тыюцхэшх юъЁєцэюёЄш $S\subset \R^2$.
═р яхЁтюь °рух яЁшьхэхэш  $E$ т $S$ чръЁр°штрхЄё  т схы√щ ЎтхЄ эх сюыхх ўхЄ√Ёхї Єюўхъ
(¤Єю Єюўъш, т ъюЄюЁ√ї ёє∙хёЄтє■Є юяюЁэ√х яЁ ь√х, ярЁрыыхы№э√х
       ъююЁфшэрЄэ√ь юё ь ш яхЁхёхър■∙шх $K$ Ёютэю т юфэющ Єюўъх).
      ┼ёыш яюёых $n$-ую °рур яЁшьхэхэш  $E$ чръЁр°хэю схы√ь ЎтхЄюь $k$ Єюўхъ,
Єю эр ёыхфє■∙хь °рух чръЁр°штрхЄё  эх сюыхх $2k$ Єюўхъ.
      ┬ ёрьюь фхых, хёыш чръЁр°штрхЄё  Єюўър $a\in E^n(S)$, Єю їюЄ  с√ эр юфэющ
       шч фтєї яЁ ь√ї, яЁюїюф ∙шї ўхЁхч $a$ ш ярЁрыыхы№э√ї ъююЁфшэрЄэ√ь юё ь,
       хёЄ№ чръЁр°хээр  Ёрэхх Єюўър.
      ┬ яЁюЄштэюь ёыєўрх ърцфр  шч ¤Єшї яЁ ь√ї т√ёхърхЄ т $E^n(S)$ эх ьхэхх
       фтєї Єюўхъ, Є.х. $a$ эх ьюцхЄ с√Є№ чръЁр°хэр эр $(n+1)$-ь °рух.
      ╚Єръ, яюёых ъюэхўэюую ўшёыр °руют сєфхЄ чръЁр°хэю ъюэхўэюх ўшёыю Єюўхъ,
       Є.х. $E^n(S)\neq \emptyset$ фы  ы■сюую $n$. QED


\begin{figure}
\begin{center}
  \begin{tabular}{@{}c@{\hskip2.5mm$\xrightarrow{\hskip5mm}$\hskip2.5mm}c@{\hskip2.5mm$\xrightarrow{\hskip5mm}$\hskip2.5mm}c@{\hskip2.5mm$\xrightarrow{\hskip5mm}$\hskip2.5mm}c@{}}
  \raisebox{-9.25mm}{\mpfile{gr}{71}}
  &
  \raisebox{-9.25mm}{\mpfile{gr}{72}}
  &
  \raisebox{0.665mm}{\mpfile{gr}{73}}
  &
  \raisebox{0.665mm}{\mpfile{gr}{74}}
  \end{tabular}
\end{center}
\caption{─юърчрЄхы№ёЄтю срчшёэющ эхтыюцшьюёЄш яхэЄюфр}\label{pic7}
\end{figure}

\smallskip
{\it ─юърчрЄхы№ёЄтю срчшёэющ эхтыюцшьюёЄш яхэЄюфр.}
╧Ёхфяюыюцшь, ўЄю яхэЄюф $T_5$ срчшёэю тыюцхэ т яыюёъюёЄ№.
╧єёЄ№ $d$ --- тхЁ°шэр яхэЄюфр.
╥ръ ъръ $E^n(T_5)=\emptyset$ фы  эхъюЄюЁюую $n$, Єю ёє∙хёЄтєхЄ ьръёшьры№эюх $k$
Єръюх, ўЄю $E^k(T_5)$ ёюфхЁцшЄ яЁюъюыюЄє■ юъЁхёЄэюёЄ№ $U$ тхЁ°шэ√ $d$ т $T_5$.
╥юуфр эр $(k+1)$-ь °рух т {\it юфэюь} шч ЁхсхЁ $A\subset T_5$
чръЁр°штрхЄё  т схы√щ ЎтхЄ эхъюЄюЁр  яюёыхфютрЄхы№эюёЄ№ Єюўхъ $a_n$,
ёїюф ∙р ё  ъ $d$.
╟эрўшЄ (яЁш эхюсїюфшьюёЄш яхЁхїюф  ъ яюфяюёыхфютрЄхы№эюёЄш т $a_n$ ш ьхэ  
эряЁртыхэшх юёш $x$), ь√ ьюцхь ёўшЄрЄ№, ўЄю юфэр шч яЁ ь√ї $x=x(a_i)$
шыш $y=y(a_i)$ эх ёюфхЁцшЄ фЁєушї Єюўхъ шч $E^k(T_5)$, ъЁюьх $a_i$.
╧юёъюы№ъє юъЁхёЄэюёЄ№ $(U-A)\cup d\cong T_4$ ёт чэр, юэр ыхцшЄ яю юфэє
ёЄюЁюэє юЄ тёхї ¤Єшї яЁ ь√ї, Є.х., т яюыєяыюёъюёЄш $\R_+\times\R$.
╧Ёш ¤Єюь $E^n(T_4)=\emptyset$, Є.х.

{\it эхъюЄюЁ√щ ъЁхёЄ $T_4\subset T_5$ срчшёэю тыюцхэ т $\R_+\times\R$ Єръ,
ўЄю $d=(0,0)$.}

╥хяхЁ№ рэрыюушўэю фюърч√трхЄё , ўЄю

{\it эхъюЄюЁ√щ ЄЁшюф $T_3\subset T_4$ срчшёэю тыюцхэ т $\R_+\times\R_+$ шыш
т $0\times \R$ Єръ, ўЄю $d=(0,0)$.}

┬ЄюЁющ ёыєўрщ эхтючьюцхэ. ╥хяхЁ№ рэрыюушўэю фюърч√трхЄё , ўЄю

{\it эхъюЄюЁ√щ фшюф $T_2\subset T_3$ срчшёэю тыюцхэ т ыєў $\R_+\times 0$ Єръ,
ўЄю $d=(0,0)$.}

╧юыєўшыш яЁюЄштюЁхўшх.
QED

\smallskip
{\it ─юърчрЄхы№ёЄтю срчшёэющ эхтыюцшьюёЄш ъЁхёЄр ё ЁрчтхЄтыхээ√ьш ъюэЎрьш.}
╧Ёхфяюыюцшь, ўЄю $C$ срчшёэю тыюцхэ т яыюёъюёЄ№.
└эрыюушўэю фюърчрЄхы№ёЄтє срчшёэющ эхтыюцшьюёЄш яхэЄюфр яюыєўрхь, ўЄю
{\it хёыш ъЁхёЄ $T_4$ срчшёэю тыюцхэ т яыюёъюёЄ№ $\R^2$, Єю юфэр шч хую тхЄтхщ
ёюфхЁцшЄ яЁ ьюышэхщэ√щ юЄЁхчюъ ё ъюэЎюь т тхЁ°шэх ъЁхёЄр, ярЁрыыхы№э√щ
юфэющ шч ъююЁфшэрЄэ√ї юёхщ.}

╥хяхЁ№ срчшёэр  эхтыюцшьюёЄ№ уЁрЇр $C$ т√ЄхърхЄ шч ёыхфє■∙хщ ыхьь√.  QED

\smallskip
{\bf ╦хььр ю ёїыюя√трэшш.}
{\it ╧єёЄ№ $K$ --- срчшёэюх яюфьэюцхёЄтю яыюёъюёЄш. 
 ╬яЁхфхышь юЄюсЁрцхэшх
$$q:\R^2\to\R^2\quad\mbox{ЇюЁьєыющ}
\quad q(x,y)=\left\{
\begin{array}{cc}
(x,y),      & x<a\\
(a,y),      & a\le x\le b \\
(x-(b-a),y),& x>b
\end{array}\right. .$$
\quad (a) $q|_{K-[a;b]\times c}$ шэ·хъЄштэю;\quad

(b) $q(K)\subset \R^2$ срчшёэю.}

\smallskip
{\it ─юърчрЄхы№ёЄтю.}
(a) ╧єёЄ№, эряЁюЄшт, фтх Єюўъш шч $K-[a;b]\times c$ ёъыхштр■Єё 
яЁш ёїыюя√трэшш $q$.
╥юуфр юэш ыхцрЄ т яюыюёх $[a;b]\times \R$ ш шьх■Є юфшэръютє■ юЁфшэрЄє $d$.
╥юуфр ¤Єш Єюўъш $(x_1,d)$, $(x_2,d)$ тьхёЄх ё $(x_1,c)$, $(x_2,c)$  ты ■Єё 
тхЁ°шэрьш яЁ ьюєуюы№эшър ёю ёЄюЁюэрьш, ярЁрыыхы№э√ьш ъююЁфшэрЄэ√ь юё ь.
▌Єю ьэюцхёЄтю эх срчшёэю т $\R^2$.
╧ЁюЄштюЁхўшх.

(b) ─юёЄрЄюўэю фюърчрЄ№, ўЄю хёыш $q(K)$ ёюфхЁцшЄ ьюыэш■
$A=\{a_1, \dots ,a_n\}$ фышэ√ $n$, Єю ш $K$ ёюфхЁцшЄ ьюыэш■ фышэ√ $n$.
┼ёыш $q^{-1}(A)$ --- ьюыэш  т $K$, Єю эєцэюх єЄтхЁцфхэшх фюърчрэю.
╚эрўх эрщфєЄё  Єюўъш $a_i,a_{i+1}$ --- тхЁ°шэ√ тхЁЄшъры№эюую чтхэр ---
Єръшх, ўЄю $p_x(a_i)=p_x(a_{i+1})$.
╥юуфр фюсртшь ъ $q^{-1}(A)$ Єюўъш $(x(q^{-1}(a_i)),c)$ ш
$(x(q^{-1}(a_{i+1})),c)$ (чфхё№ яюырурхь $q^{-1}(a,c)=(a,c)$).
╧Ёюфхырт Єръє■ юяхЁрЎш■ эхёъюы№ъю Ёрч, яюыєўшь ьюыэш■ т $K$ фышэ√ сюы№°х $n$.
QED

\bigskip
{\bf ┴рчшёэр  тыюцшьюёЄ№ т яЁюшчтхфхэшх уЁрЇют.}

{\it ─хърЁЄют√ь яЁюшчтхфхэшхь} $X\times Y$ фтєї ьэюцхёЄт $X$ ш $Y$
эрч√трхЄё  ьэюцхёЄтю тёхї ярЁ $(a,b)$ Єръшї, ўЄю $a\in X$ ш $b\in Y$.
╬яЁхфхыхэшх срчшёэюую тыюцхэш  ьюцхЄ с√Є№ юўхтшфэю юсюс∙хэю эр тыюцхэш  т
яЁюшчтюы№эюх фхърЁЄютю яЁюшчтхфхэшх $X\times Y$.
┼ёыш $X$ ш $Y$ --- уЁрЇ√, Єю ь√ ьюцхь яЁхфёЄрты Є№ ёхсх яЁюшчтхфхэшх $X\times Y$
ъръ фтєьхЁэ√щ юс·хъЄ (т эхъюЄюЁ√ї ёыєўр ї ьюцэю ёўшЄрЄ№, ўЄю ¤ЄюЄ юс·хъЄ
Ёрёяюыюцхэ т ЄЁхїьхЁэюь яЁюёЄЁрэёЄтх).
╬сючэрўшь ўхЁхч $T_n$ чтхчфє ё $n$ ыєўрьш.
═ряЁшьхЁ, фы  ЄЁшюфр $T_3$ яЁюёЄЁрэёЄтю $T_3\times I$  ты хЄё  'ъэшцъющ ё ЄЁхь 
ёЄЁрэшЎрьш', $S^1\times I$ --- ЎшышэфЁюь ш $S^1\times S^1$ --- ЄюЁюь.
╧Ёюшчтхфхэшх $T_n\times I$ эрчютхь {\it ъэшцъющ ё $n$ ёЄЁрэшЎрьш}.

\smallskip
{\bf ╟рьхўрэшх.} {\it ─ы  ы■с√ї ъюэхўэ√ї уЁрЇют $X$ ш $Y$ эрщфхЄё  ъюэхўэ√щ
уЁрЇ,  ъюЄюЁ√щ эхы№ч  срчшёэю тыюцшЄ№ т яЁюшчтхфхэшх $X\times Y$.}

─хщёЄтшЄхы№эю, юсючэрўшь ўхЁхч $k$ ьръёшьры№эє■ ёЄхяхэ№ тхЁ°шэ уЁрЇют $X$ ш
$Y$.
─юърцхь, ўЄю чтхчфр $T_{4k^2+1}$ эх тыюцшьр срчшёэю ш {\it ъєёюўэю-ышэхщэю} т
$X\times Y$.
╧Ёхфяюыюцшь, яЁюЄштэюх.
╠рыр  юъЁхёЄэюёЄ№ ЎхэЄЁр $v$ чтхчф√ т яЁюшчтхфхэшш $X\times Y$
ёюёЄюшЄ шч эх сюыхх ўхь $k^2$ яЁ ьюєуюы№эшъют $I\times J$,  ты ■∙шїё 
яЁюшчтхфхэш ьш ўрёЄхщ ЁхсхЁ $I$ ш $J$ уЁрЇют $X$ ш $Y$.
╧ю яЁшэЎшяє ─шЁшїых ёЁхфш $4k^2+1$ ЁхсхЁ чтхчф√ эрщфхЄё  яю ъЁрщэхщ ьхЁх я Є№,
эрўрыр ъюЄюЁ√ї  т√їюф Є т юфшэ яЁ ьюєуюы№эшъ.
 ╧ю¤Єюьє ёє∙хёЄтєхЄ яюфчтхчфр $T_5\subset T_{4k^2+1}$, срчшёэю
 тыюцхээр  т юфшэ шч Єръшї яЁ ьюєуюы№эшъют.
▌Єю яЁюЄштюЁхўшЄ ъЁшЄхЁш■ срчшёэющ тыюцшьюёЄш уЁрЇют.
QED

\smallskip
{\bf ╥хюЁхьр єэштхЁёры№эюёЄш.} [Ku03]
{\it ╦■сющ ъюэхўэ√щ уЁрЇ срчшёэю тыюцшь т яЁюшчтхфхэшх $X\times I$
фы  эхъюЄюЁюую сєъхЄр $X$ юъЁєцэюёЄхщ ш юЄЁхчъют.

╫шёыю юъЁєцэюёЄхщ т сєъхЄх ьюцэю тч Є№ Ёртэ√ь $E-V+C$, уфх $E,V,C$ ---
ъюышўхёЄтр тхЁ°шэ, ЁхсхЁ ш ъюьяюэхэЄ ёт чэюёЄш уЁрЇр.
┬ ўрёЄэюёЄш, ы■сюх фхЁхтю срчшёэю тыюцшью т ъэшцъє ё эхъюЄюЁ√ь ўшёыюь ёЄЁрэшЎ.}

\smallskip
╨хчєы№ЄрЄ√ ¤Єюую яєэъЄр яЁштюф Єё  схч фюърчрЄхы№ёЄтр.

╤ыхфє■∙р  ушяюЄхчр эртх эр ЄхюЁхьющ ╨юсхЁЄёюэр-╤шьюЁр ю тыюцшьюёЄш уЁрЇют т
яютхЁїэюёЄш [Sk05].

\smallskip
{\bf ├шяюЄхчр.}
{\it (a) ╤є∙хёЄтєхЄ ыш°№ ъюэхўэюх ўшёыю 'чряЁх∙хээ√ї' яюфуЁрЇют фы  срчшёэющ
тыюцшьюёЄш ъюэхўэюую уЁрЇр т фрээюх яЁюшчтхфхэшх уЁрЇют.

(b) ╤є∙хёЄтєхЄ рыуюЁшЄь яЁютхЁъш срчшёэющ тыюцшьюёЄш ъюэхўэюую уЁрЇр т
фрээюх яЁюшчтхфхэшх уЁрЇют.}

\smallskip
{\bf ╩ЁшЄхЁшщ срчшёэющ тыюцшьюёЄш уЁрЇют т ъэшцъє ё $n$ ёЄЁрэшЎрьш.} [Ku00]
─хЇхъЄюь {\it уЁрЇр $K$ эрч√трхЄё  ёєььр
$$\delta(K)=(degA_1-2)+\dots+(degA_k-2),$$
уфх $A_1,\dots ,A_k$ --- тёх тхЁ°шэ√ уЁрЇр $K$, ышсю шьх■∙шх ёЄхяхэ№ сюы№°х
ўхЄ√Ёхї, ышсю ёЄхяхэш 4, эх шьх■∙шх тшё ўшї ЁхсхЁ.
╩юэхўэ√щ уЁрЇ $K$ срчшёэю тыюцшь т $I\times T_n$ Єюуфр ш Єюы№ъю Єюуфр,
ъюуфр юэ  ты хЄё  фхЁхтюь ш

$\bullet$ ышсю $\delta(K)<n$,

$\bullet$ ышсю $\delta(K)=n$ ш $K$ ёюфхЁцшЄ тхЁ°шэє ёЄхяхэш сюы№°х ўхЄ√Ёхї,
шьх■∙є■ тшё ўхх ЁхсЁю.}

\smallskip
╚ч ¤Єюую Ёхчєы№ЄрЄр т√ЄхърхЄ яюыюцшЄхы№эюх Ёх°хэшх ушяюЄхч√ фы  яЁюшчтхфхэш 
$I\times T_n$.
┬ [Ku03] фюърчрэ√ Єръцх юсюс∙хэш  ¤Єюую Ёхчєы№ЄрЄр.

\bigskip
\newpage
\centerline{\uppercase{\bf ┴рчшёэюёЄ№ яыюёъшї ьэюцхёЄт }
\footnote{╧хЁт√щ ш юёЄры№э√х яєэъЄ√ ¤Єющ ўрёЄш эхчртшёшь√ фЁєу юЄ фЁєур.}
}


\smallskip
{\bf ╨рчЁ√тэр  срчшёэюёЄ№.}

\smallskip
{\bf 0.} ╧ЁхфёЄрт№Єх ЇєэъЎш■ $f:[(-1,-1),(1,1)]\cup[(0,0),(1,-1)]\to\R$,
$f(x,y)=xy$ т тшфх ёєьь√ $g(x)+h(y)$ фтєї ЇєэъЎшщ, ърцфр  шч ъюЄюЁ√ї чртшёшЄ
Єюы№ъю юЄ юфэющ ъююЁфшэрЄ√.

\smallskip
{\bf 1.} (a) ─ы  ы■с√ї ыш ўхЄ√Ёхї ўшёхы $f_{11},f_{12},f_{21},f_{22}$
ёє∙хёЄтє■Є Єръшх ўхЄ√Ёх ўшёыр $g_1,g_2,h_1,h_2$, ўЄю $f_{ij}=g_i+h_j$
яЁш ы■с√ї $i,j=1,2$?

(b) └эфЁхщ ═шъюырхтшў ш ┬ырфшьшЁ ╚уюЁхтшў шуЁр■Є т шуЁє '└ эє-ър, Ёрчыюцш!'.
═р °рїьрЄэющ фюёъх юЄьхўхэю эхёъюы№ъю ъыхЄюъ.
└. ═. ЁрёёЄрты хЄ ўшёыр т юЄьхўхээ√ї ъыхЄърї, ъръ їюўхЄ.
┬. ╚. ёьюЄЁшЄ эр ЁрёёЄртыхээ√х ўшёыр ш схЁхЄ 16 ўшёхы
$a_1,\dots,a_8,b_1,\dots,b_8$, Є.х. 'тхёют' ёЄюысЎют ш ёЄЁюъ, ъръ їюўхЄ.
┼ёыш ўшёыю т ърцфющ юЄьхўхээющ ъыхЄъх юърчрыюё№ Ёртэ√ь ёєььх тхёют ёЄЁюъш ш
ёЄюысЎр ¤Єющ ъыхЄъш, Єю т√шуЁры ┬. ╚., р шэрўх (Є.х. хёыш ўшёыю їюЄ  с√ т
юфэющ юЄьхўхээющ ъыхЄъх юърчрыюё№ эх Ёртэ√ь ёєььх тхёют ёЄЁюъш ш ёЄюысЎр ¤Єющ
ъыхЄъш) т√шуЁры └. ═.

─юърцшЄх, ўЄю яЁш яЁртшы№эющ шуЁх ┬. ╚. т√шуЁ√трхЄ Єюуфр ш Єюы№ъю Єюуфр, ъюуфр
эх ёє∙хёЄтєхЄ чрьъэєЄюую ьрЁ°ЁєЄр ырф№ш, эрўры№эр  ъыхЄър ш ъыхЄъш яютюЁюЄр
ъюЄюЁюую  ты ■Єё  юЄьхўхээ√ьш (эх юс чрЄхы№эю тёх юЄьхўхээ√х ъыхЄъш
чрфхщёЄтютрэ√).

\smallskip
{\it ╬яЁхфхыхэшх ьюыэшш.}
╬сючэрўшь ўхЁхч $\R^2$ яыюёъюёЄ№ ё ЇшъёшЁютрээющ ёшёЄхьющ ъююЁфшэрЄ.
╬сючэрўшь ўхЁхч $x(a)$ ш $y(a)$ ъююЁфшэрЄ√ Єюўъш $a\in\R^2$.
╧юёыхфютрЄхы№эюёЄ№ (ъюэхўэр  шыш схёъюэхўэр ) Єюўхъ яыюёъюёЄш
$\{a_1,\dots,a_n,\dots\}\subset\R^2$ эрч√трхЄё  {\it ьюыэшхщ}, хёыш фы  ърцфюую
$i$ т√яюыэхэю $a_i\neq a_{i+1}$, ш яЁш ¤Єюь $x(a_i)=x(a_{i+1})$ фы  ўхЄэ√ї
$i$ ш $y(a_i)=y(a_{i+1})$ фы  эхўхЄэ√ї $i$.
(═х юс чрЄхы№эю тёх Єюўъш ьюыэшш Ёрчышўэ√.)

╩юэхўэр  ьюыэш  $\{a_1,\dots,a_{2l+1}\}$ эрч√трхЄё  {\it чрьъэєЄющ}, хёыш
$a_1=a_{2l+1}$.

\smallskip
{\bf 2.} ╨рёёьюЄЁшь чрьъэєЄє■ ьюыэш■ $\{a_1,\dots,a_n=a_1\}$.
═рчютхь {\it Ёрчыюцхэшхь} ЁрёёЄрэютъє ўшёхы т яЁюхъЎш ї Єюўхъ ¤Єющ ьюыэшш эр
юё№ $Ox$ ш т яЁюхъЎш ї Єюўхъ ¤Єющ ьюыэшш эр юё№ $Oy$.
╠юцэю ыш Єръ ЁрёёЄртшЄ№ т Єюўърї ьюыэшш ўшёыр $f_1,\dots,f_n\in\R$ ё
$f_1=f_n$, ўЄюс√ фы  ы■сюую Ёрчыюцхэш  эхъюЄюЁюх ўшёыю $f_i$
эх с√ыю с√ Ёртэю ёєььх фтєї ўшёхы, ёЄю ∙шї т $x(a_i)$ ш т $y(a_i)$?

\smallskip
╧юфьэюцхёЄтю $K\subset\R^2$ яыюёъюёЄш эрч√трхЄё  {\it ЁрчЁ√тэю срчшёэ√ь}, хёыш
фы  ы■сющ ЇєэъЎшш $f:K\to\R$ ёє∙хёЄтє■Є Єръшх ЇєэъЎшш $g,h:\R\to\R$, ўЄю
$f(x,y)=g(x)+h(y)$ фы  ърцфющ Єюўъш $(x,y)\in K$.

\smallskip
{\bf 3.} (a) ╬ЄЁхчюъ $K=0\times[0;1]\subset\R^2$  ты хЄё  ЁрчЁ√тэю срчшёэ√ь.

(b) ╩ЁхёЄ $K=0\times[-1;1]\cup[-1;1]\times0\subset\R^2$
 ты хЄё  ЁрчЁ√тэю срчшёэ√ь.

(c) {\it ╩ЁшЄхЁшщ ЁрчЁ√тэющ срчшёэюёЄш.}
╧юфьэюцхёЄтю яыюёъюёЄш ЁрчЁ√тэю срчшёэю Єюуфр ш Єюы№ъю Єюуфр,
ъюуфр юэю эх ёюфхЁцшЄ чрьъэєЄ√ї ьюыэшщ.

\smallskip
{\bf 4.**}
└эфЁхщ ═шъюырхтшў ш ┬ырфшьшЁ ╚уюЁхтшў шуЁр■Є т 3D-шуЁє '└ эє-ър, Ёрчыюцш!'.
┬ ъєсх $n\times n\times n$, ЁрчсшЄюь эр $n^3$ хфшэшўэ√ї ъєсшъют, юЄьхўхэю
эхёъюы№ъю ъєсшъют.
└. ═. ЁрёёЄрты хЄ ўшёыр т юЄьхўхээ√ї ъєсшърї, ъръ їюўхЄ.
┬. ╚. ёьюЄЁшЄ эр ЁрёёЄртыхээ√х ўшёыр ш схЁхЄ $3n$ ўшёхы
$a_1,\dots,a_n,b_1,\dots,b_n,c_1,\dots,c_n$ --- 'тхёют' ъюыюэюъ, яЁюфюы№э√ї
ёЄЁюъ ш яюяхЁхўэ√ї ёЄЁюъ (Є.х. Ё фют, ярЁрыыхы№э√ї юёш $z$, юёш $x$ ш юёш $y$)
--- ъръ їюўхЄ.
┼ёыш ўшёыю т ърцфюь юЄьхўхээюь ъєсшъх $(i,j,k)$ (яюёЄртыхээюх └. ═.) юърчрыюё№
Ёртэ√ь ёєььх $a_i+b_j+c_k$ ЄЁхї тхёют ъюыюэъш, яЁюфюы№эющ ёЄЁюъш ш яюяхЁхўэющ
ёЄЁюъш ¤Єюую ъєсшър, Єю т√шуЁры ┬. ╚.,
р шэрўх (Є.х. хёыш ўшёыю їюЄ  с√ т юфэюь юЄьхўхээюь ъєсшъх юърчрыюё№ эх Ёртэ√ь
ёєььх ЄЁхї тхёют) т√шуЁры └. ═.

╩ръ яю эрсюЁє юЄьхўхээ√ї ъєсшъют єчэрЄ№, ъЄю т√шуЁ√трхЄ?

(▀ёэю, ўЄю {\it рыуюЁшЄь} Ёрёяючэртрэш  т√шуЁ√°эюёЄш фрээюую эрсюЁр
юЄьхўхээ√ї ъєсшъют ёє∙хёЄтєхЄ.
╞хырЄхы№эю эрщЄш яЁюёЄющ ъЁшЄхЁшщ Єшяр Єюую, ъюЄюЁ√щ шьххЄё  фы  яыюёъюую
рэрыюур ¤Єющ шуЁ√.
╚эЄхЁхёэ√ фрцх юЄтхЄ√ фы  ьрыхэ№ъшї $n$.)

\smallskip
{\bf 5.**}
(a) ╬яЁхфхышЄх ЁрчЁ√тэє■ срчшёэюёЄ№ яюфьэюцхёЄт ЄЁхїьхЁэюую яЁюёЄЁрэёЄтр.
╤ЇюЁьєышЁєщЄх ш фюърцшЄх яЁюёЄЁрэёЄтхээ√щ рэрыюу яЁштхфхээюую ъЁшЄхЁш .

(b) ╥ю цх фы  ьэюуюьхЁэюую ёыєўр .

\smallskip
{\bf ╨х°хэш  чрфрў.}

\small
\smallskip
{\bf 1.} (a) ▌Єю эхтхЁэю.
┼ёыш $f_{ij}=g_i+h_j$ фы  $i,j=1,2$, Єю
$f_{11}+f_{22}=f_{12}+f_{21}$, эю ¤Єю ёююЄэю°хэшх эх шьххЄ ьхёЄр фы 
эхъюЄюЁ√ї эрсюЁют ўшёхы $f_{ij}$.

(b) "╥юы№ъю Єюуфр" ёыхфєхЄ шч чрфрўш 2.
 ─юърцхь єЄтхЁцфхэшх "Єюуфр" шэфєъЎшхщ яю ъюышўхёЄтє юЄьхўхээ√ї ъыхЄюъ.
┼ёыш юЄьхўхэр Єюы№ъю юфэр ъыхЄър, єЄтхЁцфхэшх чрфрўш ЄЁштшры№эю.
╬сючэрўшь ўхЁхч $K$ ьэюцхёЄтю ЎхэЄЁют юЄьхўхээ√ї ъыхЄюъ.
╠эюцхёЄтю $E(K)$ юяЁхфхыхэю т ёыхфє■∙хь яєэъЄх яюёых чрфрўш 9.
╧ю єёыютш■, $K$ эх ёюфхЁцшЄ чрьъэєЄ√ї ьюыэшщ, ёыхфютрЄхы№эю $\#E(K)<\#K$.
╟эрўшЄ, яю шэфєъЄштэюьє яЁхфяюыюцхэш■ ┬.╚. ьюцхЄ т√шуЁрЄ№ эр ьэюцхёЄтх $E(K)$.
┬ёх юёЄрт°шхё  ъыхЄъш  ты ■Єё  хфшэёЄтхээ√ьш юЄьхўхээ√ьш т ётюхщ ёЄЁюъх шыш
т ётюхь ёЄюысЎх.
╤ыхфютрЄхы№эю, ┬.╚. ёьюцхЄ т√сЁрЄ№ ш юёЄрт°шхё  тхёр фы  $K$.

\smallskip
{\bf 2.} ─р. ┼ёыш ърцфюх шч ўшёхы $f_i$ яЁхфёЄртшью т тшфх ёєьь√
фтєї ўшёхы, Ёрёяюыюцхээ√ї т Єюўърї $x(a_i)$ ш $y(a_i)$, Єю
 $f_1-f_2+f_3- \dots -f_{n-1}=0$,
эю ьюцэю ыхуъю яюфюсЁрЄ№ эрсюЁ ўшёхы $f_i$, фы  ъюЄюЁюую ¤Єю эхтхЁэю.

\smallskip
{\bf 3.} (a)  ╧юыюцшь $h(y)=f(0,y)$ ш $g(x)=0$.

(b) ╧юыюцшь $g(x)=f(x,0)$ ш $h(y)=f(0,y)-f(0,0)$.

(c) "╥юы№ъю Єюуфр" ёыхфєхЄ шч чрфрўш 2.
─юърцхь єЄтхЁцфхэшх "Єюуфр".
╨рёёьюЄЁшь яЁюшчтюы№эє■ ЇєэъЎш■ $f:K\to \R$ ш яюёЄЁюшь яю эхщ ЇєэъЎшш
$g$ ш $h$ Єръшх, ўЄю $f(x,y)=g(x)+h(y)$.
═рчют╕ь фтх Єюўъш $a,b\in K$
{\it ¤ътштрыхэЄэ√ьш}, хёыш ёє∙хёЄтєхЄ ьюыэш  $\{a=a_1,\dots,a_n=b\}\subset K$.
┬юч№ь╕ь юфшэ шч ъырёёют ¤ътштрыхэЄэюёЄш $K_1\subset K$ ш юяЁхфхышь
ЇєэъЎшш $g:x(K_1)\to\R$ ш $h:y(K_1)\to\R$ ёыхфє■∙шь юсЁрчюь. ╟рЇшъёшЁєхь
яЁюшчтюы№эє■ Єюўъє $a_1\in K_1$. ╧юыюцшь $g(x(a_1))=f(a_1)$ ш $h(y(a_1))=0$.
┼ёыш $\{a_1,a_2,\dots,a_{2l}\}$ --- ьюыэш  шч Єюўхъ ьэюцхёЄтр $K$, Єю
яюыюцшь
$$h(y(a_{2l})):=f(a_{2l})-f(a_{2l-1})+\dots -f(a_1)\quad\mbox{ш}
\quad g(x(a_{2l})):=f(a_{2l-1})-f(a_{2l-2})+\dots+f(a_1).$$
┼ёыш $\{a_1,a_2,\dots,a_{2l+1}\}$ --- ьюыэш  шч Єюўхъ ьэюцхёЄтр $K$, Єю
яюыюцшь
$$g(x(a_{2l+1}))=f(a_{2l+1})-f(a_{2l})+\dots+f(a_1)$$
(чэрўхэшх $h(y(a_{2l+1}))$ єцх юяЁхфхыхэю).
╤фхырхь ¤Єю яюёЄЁюхэшх фы  тёхї ъырёёют ¤ътштрыхэЄэюёЄш юфэютЁхьхээю.
─ы  тёхї цх яЁюўшї Єюўхъ яюыюцшь $g(x)=0$ ш $h(y)=0$.

\normalsize

\bigskip
{\bf ═хяЁхЁ√тэр  срчшёэюёЄ№.}

\smallskip
╧юфьэюцхёЄтю $K\subset\R^2$ эрч√трхЄё  {\it (эхяЁхЁ√тэю) срчшёэ√ь}, хёыш
фы  ы■сющ эхяЁхЁ√тэющ ЇєэъЎшш $f:K\to\R$ ёє∙хёЄтє■Є Єръшх эхяЁхЁ√тэ√х ЇєэъЎшш
$g,h:\R\to\R$, ўЄю $f(x,y)=g(x)+h(y)$ фы  ърцфющ Єюўъш $(x,y)\in K$.
(╬яЁхфхыхэшх эхяЁхЁ√тэющ ЇєэъЎшш эряюьэхэю т эрўрых ўрёЄш 1.)
╤ыютю 'эхяЁхЁ√тэю' фрыхх юяєёърхЄё .

\smallskip
{\bf ╧Ёюсыхьр └Ёэюы№фр.} {\it ╩ръшх яюфьэюцхёЄтр яыюёъюёЄш  ты ■Єё  срчшёэ√ьш?}

\smallskip
╫Єюс√ яюфющЄш ъ юЄтхЄє, ЁрёёьюЄЁшь эхёъюы№ъю яЁшьхЁют.

\smallskip
{\bf 6.} (a) ╟рьъэєЄр  ьюыэш  эх срчшёэр.

(b) ╬ЄЁхчюъ $K=0\times[0;1]\subset\R^2$  ты хЄё  срчшёэ√ь.

(c) ╩ЁхёЄ $K=0\times[-1;1]\cup[-1;1]\times0\subset\R^2$  ты хЄё  срчшёэ√ь.

\smallskip
{\bf 7.}
(a) ┼ёыш яюфьэюцхёЄтю яыюёъюёЄш срчшёэю, Єю юэю ЁрчЁ√тэю срчшёэю.
(╬яЁхфхыхэшх ш ъЁшЄхЁшщ ЁрчЁ√тэющ срчшёэюёЄш ёь. т яЁхф√фє∙хь яєэъЄх.)

(b) {\it ╧юяюыэхээющ ьюыэшхщ} эрч√трхЄё  юс·хфшэхэшх Єюўъш $a_0\in\R^2$ ё
схёъюэхўэющ ьюыэшхщ $\{a_1,\dots,a_n,\dots\}\subset\R^2$ шч Ёрчышўэ√ї Єюўхъ,
{\it ёїюф ∙хщё } ъ Єюўъх $a_0$ (Є.х. фы  ы■сюую $\eps>0$ эрщфхЄё  Єръюх
эрЄєЁры№эюх $N$, ўЄю фы  ы■сюую $i>N$ т√яюыэхэю $|a_i,a_0|<\eps$).
─юърцшЄх, ўЄю эшърър  яюяюыэхээр  ьюыэш  эх  ты хЄё  срчшёэющ.
(╟рьхЄшь, ўЄю юэр  ты хЄё  ЁрчЁ√тэю срчшёэющ).

(c) ╫хЁхч $[a,b]$ юсючэрўшь юЄЁхчюъ, ёюхфшэ ■∙шщ Єюўъш $a$ ш $b$.
─юърцшЄх, ўЄю ъЁхёЄ 
\linebreak
$[(-1,-2),(1,2)]\cup[(-1,1),(1,-1)]$ эх  ты хЄё 
срчшёэ√ь.

(d) ╧єёЄ№ $m_{i,j}=2-3\cdot2^{-i}+j\cdot2^{-2i}$.
╨рёёьюЄЁшь ьэюцхёЄтю, ёюёЄю ∙хх шч Єюўхъ $(m_{i,2l},m_{i,2l})$ ш Єюўхъ
$(m_{i,2l},m_{i,2l-2})$, уфх $i=1,2,\dots$ ш $l=1,2,3,\dots,2^{i-1}$.
─юърцшЄх, ўЄю ¤Єю яюфьэюцхёЄтю яыюёъюёЄш эх ёюфхЁцшЄ схёъюэхўэющ ьюыэшш,
эю ёюфхЁцшЄ ёъюы№ єуюфэю фышээ√х ьюыэшш.

(e) ╬с·хфшэхэшх ьэюцхёЄтр шч яЁхф√фє∙хую яєэъЄр ё Єюўъющ $(2,2)$ эх срчшёэю.

\smallskip
╧юёыхфютрЄхы№эюёЄ№ Єюўхъ $a_i$ яыюёъюёЄш эрч√трхЄё  {\it ёїюф ∙хщё  ъ Єюўъх
$a$}, хёыш фы  ы■сюую $\eps>0$ эрщфхЄё  Єръюх Ўхыюх $N$, ўЄю фы  ы■сюую $i>N$
т√яюыэхэю $|a,a_i|<\eps$.

╧юфьэюцхёЄтю $K\subset\R^2$ яыюёъюёЄш эрч√трхЄё  {\it чрьъэєЄ√ь},
хёыш фы  ы■сющ схёъюэхўэющ яюёыхфютрЄхы№эюёЄш Єюўхъ $a_i\in K$,
ёїюф ∙хщё  ъ Єюўъх $a$, т√яюыэхэю $a\in K$.

\smallskip
{\bf 8.} ╧юфьэюцхёЄтю $K\subset\R^2$ яыюёъюёЄш  ты хЄё 
чрьъэєЄ√ь Єюуфр ш Єюы№ъю Єюуфр, ъюуфр фы  ы■сющ Єюўъш $a\not\in K$
эрщфхЄё  Єръюх $\eps>0$, ўЄю ы■ср  Єюўъa яыюёъюёЄш c ЁрёёЄю эшхь
ьхэхх $\eps$ фю $a$ эх яЁшэрфыхцшЄ $K$.


\smallskip
{\bf ╩ЁшЄхЁшщ срчшёэюёЄш.}
{\it ╟рьъэєЄюх юуЁрэшўхээюх яюфьэюцхёЄтю яыюёъюёЄш срчшёэю Єюуфр ш
Єюы№ъю Єюуфр, ъюуфр юэю эх ёюфхЁцшЄ ёъюы№ єуюфэю фышээ√ї ьюыэшщ} [St89].

\smallskip
╧Ёштхфхь чфхё№ чрьхўрэш  ш яхЁхЇюЁьєышЁютъє (шёяюы№чєхьє■ т фюърчрЄхы№ёЄтх).
╤рью фюърчрЄхы№ёЄтю яЁштюфшЄё  т ёыхфє■∙хь яєэъЄх.

\smallskip
{\bf 9.}
(a) ╙ёыютшх чрьъэєЄюёЄш т ъЁшЄхЁшш фхщёЄтшЄхы№эю эхюсїюфшью (Є.х. хёыш т
ЇюЁьєышЁютъх ЄхюЁхь√ юяєёЄшЄ№ ¤Єю єёыютшх, Єю яюыєўшЄё  эхтхЁэюх єЄтхЁцфхэшх).

(b) ╙ёыютшх юуЁрэшўхээюёЄш т ъЁшЄхЁшш фхщёЄтшЄхы№эю эхюсїюфшью (Є.х. хёыш т
ЇюЁьєышЁютъх ЄхюЁхь√ юяєёЄшЄ№ ¤Єю єёыютшх, Єю яюыєўшЄё  эхтхЁэюх єЄтхЁцфхэшх).

(ё)** ═рщфшЄх ъЁшЄхЁшщ срчшёэюcЄш фы  чрьъэєЄ√ї (эю эхюуЁрэшўхээ√ї)
яюфьэюцхёЄт яыюёъюёЄш.


\smallskip
╧єёЄ№ $K$ --- яюфьэюцхёЄтю яыюёъюёЄш $\R^2$.
─ы  ърцфющ Єюўъш $v\in K$ эрЁшёєхь фтх яЁ ь√х, яЁюїюф ∙шх ўхЁхч $v$
ярЁрыыхы№эю ъююЁфшэрЄэ√ь юё ь.  ┼ёыш їюЄ  с√ юфэр шч ¤Єшї фтєї яЁ ь√ї
яхЁхёхърхЄ $K$ Єюы№ъю т Єюўъх $v$, Єю яюъЁрёшь $v$ т схы√щ ЎтхЄ.  ╬сючэрўшь
ўхЁхч $E(K)$ ьэюцхёЄтю тёхї Єюўхъ $K$, эх  ты ■∙шїё  схы√ьш:
$$E(K)=\{v\in K:\ |K\cap(x=x(v))|\ge2\mbox{ ш }|K\cap(y=y(v))|\ge2\}.$$
═ряЁшьхЁ, Ёшё. 6b яюыєўрхЄё  шч Ёшё. 6a юяхЁрЎшхщ $E$. 
╧єёЄ№ $E^2(K)=E(E(K))$, $E^3(K)=E(E(E(K)))$ ш Є.ф.

\smallskip
{\bf 10.}
╧юфьэюцхёЄтю $K\subset\R^2$ эх ёюфхЁцшЄ ёъюы№ єуюфэю фышээ√ї ьюыэшщ Єюуфр ш
Єюы№ъю Єюуфр, ъюуфр  $E^n(K)=\emptyset$ фы  эхъюЄюЁюую $n$.

\smallskip
{\bf 12.} (a)* ─юърцшЄх ¤ыхьхэЄрЁэю (Є.х. схч шёяюы№чютрэш  юяшёрэш 
яЁюёЄЁрэёЄтр $C^*(K)$ т ЄхЁьшэрї ьхЁ, ёь. ёыхфє■∙шщ яєэъЄ), 
ўЄю хёыш $K\subset\R^2$
чрьъэєЄю ш юуЁрэшўхэю, яЁшўхь $E(K)=\emptyset$, Єю $K$ срчшёэю [Mi09].

╙ърчрэшх. ╧юыєўшЄх ёэрўрыр Ёрчыюцхэшх $f(x,y)=g(x)+h(y)$ фы 
{\it ъєёюўэю-ышэхщэ√ї} ЇєэъЎшщ $f$, яЁшўхь $|g|+|h|<5|f|$.

(b)** ─юърцшЄх ¤ыхьхэЄрЁэю ўрёЄ№ `Єюуфр' ъЁшЄхЁш  срчшёэюёЄш.

╙ърчрэшх. ╥ю цх, $|g|+|h|<C_n|f|$, уфх $C_n$ чртшёшЄ Єюы№ъю юЄ Єюую $n$, фы 
ъюЄюЁюую $E^n(K)=\emptyset$.


\smallskip
{\bf 11.}
┴рчшёэюёЄ№ яюфьэюцхёЄт ЄЁхїьхЁэюую яЁюёЄЁрэёЄтр юяЁхфхыхэр т√°х яхЁхф ыхььющ
└Ёэюы№фр ю фхЁхт№ ї.

(a)
─юърцшЄх, ўЄю хц $0\times0\times[-1;1]\cup0\times[-1;1]\times0\cup
[-1;1]\times0\times0\subset\R^3$  ты хЄё  срчшёэ√ь.

(b) ╧юфьэюцхёЄтю яЁюёЄЁрэёЄтр $\R^3$, ёюёЄю ∙хх шч ўхЄ√Ёхї Єюўхъ
$(0,0,0)$; $(1,1,0)$; $(0,1,1)$; $(1,0,1)$, срчшёэю.
(═ю $E^n(K)\neq\emptyset$ фы   ы■сюую $n$, ёь. эшцх.)

(c)* ─ы  $K\subset\R^3$ рэрыюушўэю юяЁхфхышь $E(K)$, шёяюы№чє  тьхёЄю яЁ ь√ї
яыюёъюёЄш, яхЁяхэфшъєы Ёэ√х  юё ь ъююЁфшэрЄ:
$$E(K):=\{v\in K:\ |K\cap(x=x(v))|\ge2,
\ |K\cap(y=y(v))\ge2|\mbox{ ш }|K\cap(z=z(v))|\ge2\}.$$
─юърцшЄх, ўЄю хёыш $K$ чрьъэєЄю, юуЁрэшўхэю ш $E^n(K)=\emptyset$ фы  эхъюЄюЁюую 
$n$, Єю $K$ срчшёэю [St89, Lemma 23.ii].

(d)* ═шъръюх яюфьэюцхёЄтю яЁюёЄЁрэёЄтр $\R^3$ (шыш фрцх $\R^4$), уюьхюьюЁЇэюх
фтєьхЁэюьє фшёъє, эх  ты хЄё  срчшёэ√ь.
╙ърчрэшх: юяЁхфхыхэшх ьэюуюьхЁэющ ьюыэшш ёь. т [St89, 6.12, p.39].

\smallskip
{\bf ╨х°хэш  чрфрў.}

\small
\smallskip
{\bf 6.}
(a) ┼ёыш с√ ьюыэш  $A=\{a_1,\dots,a_{2l+1}=a_1\}$ с√ыр срчшёэющ, Єю
$f(a_1)-f(a_2)+\dots +f(a_{2l-1})-f(a_{2l})=0$, эю ыхуъю яюфюсЁрЄ№
ЇєэъЎш■ $f$, фы  ъюЄюЁющ ¤Єю эх т√яюыэхэю. ╤ЁртэшЄх ё чрфрўхщ 2.

(b),(c) └эрыюушўэю чрфрўрь 3a, 3b.

\smallskip
{\bf 7.} (a) ┼ёыш ьэюцхёЄтю эх  ты хЄё  ЁрчЁ√тэю срчшёэ√ь, Єю яю ъЁшЄхЁш■
ЁрчЁ√тэющ срчшёэюёЄш шч яЁхф√фє∙хую яєэъЄр юэю ёюфхЁцшЄ чрьъэєЄє■ ьюыэш■.
╥юуфр єЄтхЁцфхэшх чрфрўш ёыхфєхЄ шч 6a, Єръ ъръ ЇєэъЎш  $f$ ьюцхЄ
с√Є№ яЁюфюыцхэр ё чрьъэєЄющ ьюыэшш эр тё╕ ьэюцхёЄтю.

(b) ╨рёёьюЄЁшь ЇєэъЎш■ $f$, фы  ъюЄюЁющ $f(a_i)=\frac{(-1)^i}i$.
╧Ёхфяюыюцшь, ўЄю $f(x,y)=g(x)+h(y)$ фы  эхъюЄюЁ√ї эхяЁхЁ√тэ√ї $g$ ш $h$,
Єюуфр
$$f(a_1)-f(a_2)+f(a_3)-f(a_4)+ \dots-f(a_{2l})=h(y(a_1))-h(y(a_{2l})).$$
╥ръ ъръ $\lim_{l\to\infty}h(y_{2l})$ ёє∙хёЄтєхЄ ш Ёртхэ
$h(y(a_0))$, Єю Ё ф $\sum_{i=1}^{2l} (-1)^i f(a_i)$ ёїюфшЄё 
яЁш $l\to\infty$. ═ю ¤Єю яЁюЄштюЁхўшЄ ЁрёїюфшьюёЄш урЁьюэшўхёъюую Ё фр.

(c) ╩ЁхёЄ ёюфхЁцшЄ чрьъэєЄє■ ьюыэш■
$$a_{4k+1}=(\frac{-1}{4^k},\frac{1}{4^k}),
\ a_{4k+2}=(\frac{1}{2\cdot4^k},\frac{1}{4^k}),
\ a_{4k+3}=(\frac{1}{2\cdot4^k},\frac{-1}{2\cdot4^k}),
\ a_{4k+4}=(\frac{-1}{4^{k+1}},\frac{-1}{2\cdot4^k})$$
╬яЁхфхышь ЇєэъЎш■ $f$ эр ¤Єющ ьюыэшш, шёяюы№чє  чрфрўє 7(b), ш яЁюфюыцшь
х╕ ъєёюўэю-ышэхщэю эр тхё№ ъЁхёЄ. ═х ёє∙хёЄтєхЄ Єръшї ЇєэъЎшщ $g$ ш $f$, ўЄю
$f(x,y)=g(x)+h(y)$.

(d) ─ы  ы■сюую $i$ Єюўъш $(m_{i,2l},m_{i,2l})_{l=1}^{2^{i-1}}$ ш
$(m_{i,2l},m_{i,2l-2})_{l=1}^{2^{i-1}}$
юсЁрчє■Є ьюыэш■ шч $2^i$ ¤ыхьхэЄют.

(e) ╬яЁхфхышь ЇєэъЎш■ $f(x,y)$ ёююЄэю°хэш ьш
$$f((m_{i,2l},m_{i,2l})):=\frac 1{2^i} \quad \mbox{ ш } \quad
f((m_{i,2l},m_{i,2l-2})):=-\frac 1{2^i}$$
╧Ёхфяюыюцшь, ўЄю $f(x,y)=g(x)+h(y)$
фы  эхъюЄюЁ√ї эхяЁхЁ√тэ√ї $g(x)$ ш $h(y)$.
╥хяхЁ№ фы  ърцфюую $i$, шёяюы№чє  ьюыэшш $(m_{i,2l},m_{i,2l})$ ш
 $(m_{i,2l},m_{i,2l-2})$, уфх $l=1,2,3,\dots 2^{i-1}$, яюыєўрхь
$h(2-\frac 3{2^i})-h(2-\frac 2{2^i})=1$.
 ▌Єю яЁюЄштюЁхўшЄ эхяЁхЁ√тэюёЄш $h$ т Єюўъх $y=2$.

\smallskip
{\bf 8.} ─юърцхь єЄтхЁцфхэшх "Єюы№ъю Єюуфр".
╧єёЄ№ $K$ --- чрьъэєЄюх яюфьэюцхёЄтю яыюёъюёЄш.
╧Ёхфяюыюцшь, ўЄю фы  эхъюЄюЁющ Єюўъш $a=(x,y)\not\in K$
ш фы  яЁюшчтюы№эюую $\eps=\frac 1n>0$
ёє∙хёЄтєхЄ їюЄ  с√ юфэр Єюўър $a_n\in K$, фы  ъюЄюЁющ  $|a,a_n| \le \frac 1n$.
═ю Єюуфр яюёыхфютрЄхы№эюёЄ№ Єюўхъ $a_n\in K$ ёїюфшЄё  ъ Єюўъх $a$, яю¤Єюьє $a\in K$.
╧ЁюЄштюЁхўшх.

╥хяхЁ№ фюърцхь єЄтхЁцфхэшх "Єюуфр".
╧єёЄ№ эхъюЄюЁр  яюёыхфютрЄхы№эюёЄ№ $a_n$ ёїюфшЄё  ъ Єюўъх $a$, эх ыхцр∙хщ
т ьэюцхёЄтх $K$. ╧ю єёыютш■ ёє∙хёЄтєхЄ $\eps>0$ Єръюх, ўЄю фы  ы■сющ
Єюўъш $a_n \in K$ ЁрёёЄю эшх $|a,a_n|>\eps$. ═ю ¤Єю яЁюЄштюЁхўшЄ ёїюфшьюёЄш
яюёыхфютрЄхы№эюёЄш.

\smallskip
{\bf 9.}
(a) ╦■ср  схёъюэхўэр  ьюыэш  $A$, эх ёюфхЁцр∙р 
чрьъэєЄ√ї ьюыэшщ ш ёїюф ∙р ё  ъ Єюўъх $a\not\in A$,  ты хЄё  срчшёэющ.
▌Єю ёыхфєхЄ шч Єюую, ўЄю ы■ср  ЇєэъЎш , юяЁхфхы╕ээр  эр $A$, эхяЁхЁ√тэр.

(b) ╩юэЄЁяЁшьхЁюь  ты хЄё  ьэюцхёЄтю
$\{(k,k)\}_{k=1}^\infty\cup\{(k,k-1)\}_{k=1}^\infty$ Єюўхъ яыюёъюёЄш.

\smallskip
{\bf 10.}
─юърцхь ўрёЄ№ 'Єюы№ъю Єюуфр'.
╧Ёхфяюыюцшь, ўЄю $E^n(K)\neq\emptyset$ фы  тёхї $n$.
─ы  ърцфюую $n$ ЁрёёьюЄЁшь Єюўъє $a_0\in E^n(K)$. ┬√схЁхь Єюўъш $a_{-1},a_1\in
E^{n-1}(K)$ Єръшх, ўЄю $x(a_{-1})=x(a_0)$ ш $y(a_1)=y(a_0)$.
╥хяхЁ№ ьюцэю т√сЁрЄ№ Єюўъш $a_{-2},a_2\in E^{n-2}(K)$, фы  ъюЄюЁ√ї
$\{a_{-2},a_{-1},a_0,a_1,a_2\}$ --- ьюыэш . └эрыюушўэю ьюцэю ёъюэёЄЁєшЁютрЄ№
ьюыэш■ шч $2n+1$ Єюўхъ, ыхцр∙є■ Ўхышъюь т ьэюцхёЄтх $K$. ╫Єю ш ЄЁхсютрыюё№ фюърчрЄ№.

─юърцхь ўрёЄ№ 'Єюуфр'.
╧єёЄ№ ьэюцхёЄтю $K$ ёюфхЁцшЄ ьюыэш■ шч $2n+1$ Єюўъш
 $\{a_{-n},\dots,a_0,\dots,a_{n}\}$. ╥юуфр т ьэюцхёЄтх $E(K)$ ёюфхЁцшЄё 
ьюыэш  шч $2n-1$ Єюўъш
 $\{a_{-n+1},\dots,a_{n-1}\}$. ╧Ёюфюыцр , яюыєўшь, ўЄю $a_0\in E^n(K)$.
╤ыхфютрЄхы№эю, хёыш $E^n(K)=\emptyset$, Єю $K$ эх ёюфхЁцшЄ
ьюыэшш шч $2n+1$ Єюўхъ.


\smallskip
{\bf 11.}
(a) ─ы  яЁюшчтюы№эющ ЇєэъЎшш $f:K \to \R$ эр хцх $K$ юяЁхфхышь
 $g(x):=f(x,0,0)$, $h(y):=f(0,y,0)-f(0,0,0)$ ш $l(z):=f(0,0,z)-f(0,0,0)$.

(b) ╧юыюцшь $g(0)=f(0,0,0)$, $h(0)=0$, $l(0)=0$,
$$2g(1)=f(0,0,0)+f(1,1,0)+f(1,0,1)-f(0,1,1),$$
$$2h(1)=-f(0,0,0)+f(1,1,0)-f(1,0,1)+f(0,1,1))\quad\mbox{ш}$$
$$2l(1)=-f(0,0,0)-f(1,1,0)+f(1,0,1)+f(0,1,1).$$

\normalsize

\bigskip
{\bf ─юърчрЄхы№ёЄтю ъЁшЄхЁш  срчшёэюёЄш.}

╧єёЄ№ $K$ --- яЁюшчтюы№эюх чрьъэєЄюх юуЁрэшўхээюх яюфьэюцхёЄтю яыюёъюёЄш.
╚чтхёЄэю, ўЄю Єюуфр ы■ср  эхяЁхЁ√тэр  ЇєэъЎш  $f:K\to \R$ юуЁрэшўхэр.
╘єэъЎш  $f:K\to \R$ эрч√трхЄё  {\it юуЁрэшўхээющ}, хёыш эрщфхЄё  ўшёыю $M$
Єръюх, ўЄю $|f(x)|<M$  фы  ы■сющ Єюўъш $x\in K$.
─ы  юуЁрэшўхээющ ЇєэъЎшш $G:K \to \R$ яюыюцшь $|G|:=\sup_{x\in K}|G(x)|$.

\smallskip
{\it ═рўрыю фюърчрЄхы№ёЄтр ўрёЄш 'Єюы№ъю Єюуфр' ъЁшЄхЁш  срчшёэюёЄш.}
╧Ёхфяюыюцшь, эряЁюЄшт, ўЄю $K$ ёюфхЁцшЄ ёъюы№ єуюфэю фышээ√х ьюыэшш ш срчшёэю.
┬√сшЁр  яюфяюёыхфютрЄхы№эюёЄш, ьюцэю фюсшЄ№ё  Єюую, ўЄюс√ т ърцфющ ьюыэшш
Єюўъш яюярЁэю Ёрчышўэ√.
╧ю¤Єюьє сєфхь ёўшЄрЄ№, ўЄю ¤Єю т√яюыэхэю.
╥юуфр фы  ы■сюую $n\ge4$ ёє∙хёЄтєхЄ ьюыэш  $\{a^n_1,\dots,a^n_{2n+5}\}$ шч
$(2n+5)$-ш Ёрчышўэ√ї Єюўхъ ьэюцхёЄтр $K$.

╤є∙хёЄтєхЄ эхяЁхЁ√тэр  ЇєэъЎш 
$$f_n:K\to\R\quad\mbox{Єрър , ўЄю}\quad f_n(a^n_i)=(-1)^i\quad\mbox{ш}\quad
|f_n(x)|\le1\quad\mbox{фы  ы■сюую}\quad x\in K.$$
(─хщёЄтшЄхы№эю, яюёЄЁюшь ёэрўрыр эхяЁхЁ√тэє■ ЇєэъЎш■ $f:\R^2\to \R$,
єфютыхЄтюЁ ■∙є■ ¤Єшь єёыютш ь.
╬сючэрўшь $s=\min_{i<j}|a_i,a_j|$.
╨рёёьюЄЁшь $n$ фшёъют ё ЎхэЄЁрьш т Єюўърї $a_i$ ш Ёрфшєёрьш $\frac s3$.
┬эх ¤Єшї фшёъют яюыюцшь $f=0$. ┬эєЄЁш $i$-ую фшёър ёфхырхь
$f$ ышэхщэющ ЇєэъЎшхщ юЄ Ёрфшєёр, Ёртэющ $(-1)^i$ т ЎхэЄЁх $a_i$ ш эєы■
эр уЁрэшЎх.
╥хяхЁ№ юуЁрэшўшь яюёЄЁюхээє■ ЇєэъЎш■
эр $K\subset \R^2$ ш яюыєўшь ЄЁхсєхьє■ эхяЁхЁ√тэє■ ЇєэъЎш■ $K\to \R$.)

╬яЁхфхышь яю шэфєъЎшш яюёыхфютрЄхы№эюёЄ№ ўшёхы $s_n$ ш ЇєэъЎшщ $F_n:K\to \R$.
╧юыюцшь $s_0=1$ ш $F_0=0$. ╧Ёхфяюыюцшь, ўЄю $s_{n-1}$ ш $F_{n-1}$ єцх
юяЁхфхыхэ√.
┬юч№ьхь ЇєэъЎшш $G_{n-1},H_{n-1}:\R\to\R$, фы  ъюЄюЁ√ї
$F_{n-1}(x,y)=G_{n-1}(x)+H_{n-1}(y)$ (хёыш Єръшї ЇєэъЎшщ эхЄ, Єю тёх фюърчрэю).
┴хЁхь
$$s_n>s_{n-1}!\cdot(|G_{n-1}|+n)\quad\mbox{ш}
\quad F_n=F_{n-1}+\frac{f_{s_n}}{s_{n-1}!}$$
─юёЄрЄюўэю фюърчрЄ№, ўЄю ЇєэъЎш 
$$\displaystyle F=\sum\limits_{n=1}^\infty\frac{f_{s_n}}{s_{n-1}!}=
\lim\limits_{n\to\infty}F_n$$
эх яЁхфёЄртшьр т тшфх $G(x)+H(y)$.

╧Ёхфяюыюцшь, ўЄю, эряЁюЄшт $F(x,y)=G(x)+H(y)$ фы  эхъюЄюЁ√ї $G$ ш $H$.
─ы  яюыєўхэш  яЁюЄштюЁхўш  фюёЄрЄюўэю фюърчрЄ№, ўЄю $|G|>n$ фы  ърцфюую $n$.
└ фы  ¤Єюую фюёЄрЄюўэю яюърчрЄ№, ўЄю $s_{n-1}!|G-G_{n-1}|>s_n$: Єюуфр сєфхЄ
$$|G|+|G_{n-1}|\ge|G-G_{n-1}|>\frac{s_n}{s_{n-1}!}>|G_{n-1}|+n.$$


{\bf ╦хььр.}
{\it ╧єёЄ№ $m\ge4$,

$\bullet$ $K=\{a_1,\dots,a_{2m+5}\}$ --- ьюыэш  шч $2m+5$ Ёрчышўэ√ї Єюўхъ эр
яыюёъюёЄш,

$\bullet$ $f(a_1),\dots,f(a_{2m+5})$ --- ўшёыр, фы  ъюЄюЁ√ї
$|(-1)^i-f(a_i)|<1/m$ ш

$\bullet$ $g(x(a_i)),h(y(a_i))$, $i=1,\dots,2m+5$, --- Єръшх ўшёыр, ўЄю
$f(a_i)=g(x(a_i))+h(y(a_i))$ фы  ы■сюую $i$ (яЁш ¤Єюь хёыш $x(a_i)=x(a_j)$, Єю
$g(x(a_i))=g(x(a_j))$, ш рэрыюушўэю фы  $y$ ш $h$).

╥юуфр $\max_i|g(x(a_i))|>m$.}


\smallskip
{\it ─юърчрЄхы№ёЄтю.}
╠юцэю ёўшЄрЄ№, ўЄю яЁ ьр  $a_1a_2$ ярЁрыыхы№эр юёш $Ox$ (хёыш ¤Єю эх Єръ,
Єю єтхышўшь тёх шэфхъё√ эр 1 т яюёыхфє■∙шї ЇюЁьєырї).
╚ьххь
$$|(f(a_1)-f(a_2)+f(a_3)-f(a_4)+\dots-f(a_{2m+4}))-(2m+4)|\le\frac{2m+4}m\le3.$$
▌Єю ючэрўрхЄ, ўЄю $|g(x(a_1))-g(x(a_{2m+4}))|\ge(2m+4)-3>2m$.
╬Єё■фр ёыхфєхЄ ЄЁхсєхьюх эхЁртхэёЄтю.
QED

\smallskip
{\it ╬ъюэўрэшх фюърчрЄхы№ёЄтр ўрёЄш 'Єюы№ъю Єюуфр' ъЁшЄхЁш  эхяЁхЁ√тэющ
срчшёэюёЄш.}
╚ьххь:
$$F-F_n=F-F_{n-1}-\frac{f_{s_n}}{s_{n-1}!}=
\frac{s_{n-1}!(F-F_{n-1})-f_{s_n}}{s_{n-1}!}.$$
╧Ёшьхэшь ыхььє ъ
$$m=s_n,\quad a_i=a_i^{s_n},\quad f=s_{n-1}!(F-F_{n-1}), \quad
g=s_{n-1}!(G-G_{n-1}),\quad h=s_{n-1}!(H-H_{n-1}).$$
▌Єю тючьюцэю, Єръ ъръ $f(x,y)=g(x)+h(y)$ ш (Єръ ъръ $s_n-1>s_{n-1}$ яЁш $n>2$)
$$|f-f_{s_n}|=s_{n-1}!|F-F_n|<\frac1{(s_n-1)\cdot s_n}\sum\limits_{k=0}^\infty
\frac1{(s_n+1)\cdot\dots\cdot s_{n+k}}<
\frac1{(s_n-1)\cdot s_n}\sum\limits_{k=0}^\infty\frac1{2^k}<\frac1{s_n}.$$
╧ю ыхььх яюыєўшь $s_{n-1}!|G-G_{n-1}|>s_n$.
QED

\smallskip
{\it ─юърчрЄхы№ёЄтю ъЁшЄхЁш  срчшёэюёЄш} [St89, \S2, ╦хььр 23.ii].
\footnote{▌Єю фюърчрЄхы№ёЄтю эх¤ыхьхэЄрЁэю, ЇюЁьры№эю эх шёяюы№чєхЄё  т
фры№эхщ°хь ш ьюцхЄ с√Є№ юяє∙хэю ўшЄрЄхыхь.
╬фэръю ь√ яЁштюфшь хую, яюёъюы№ъє эр°х шчыюцхэшх ъюЁюўх ш  ёэхх фрээюую т
[St89].}
╬эю юёэютрэю эр яхЁхЇюЁьєышЁютъх ётющёЄтр срчшёэюёЄш т ЄхЁьшэрї
{\it юуЁрэшўхээ√ї ышэхщэ√ї юяхЁрЄюЁют} т {\it срэрїют√ї яЁюёЄЁрэёЄтрї ЇєэъЎшщ}.
╬сючэрўшь ўхЁхч $C(X)$ яЁюёЄЁрэёЄтю
эхяЁхЁ√тэ√ї ЇєэъЎшщ эр $X$ ё эюЁьющ $|f|=\sup\limits\{|f(x)|\ :\ x\in X\}$.
┬ ¤Єюь фюърчрЄхы№ёЄтх юсючэрўшь ўхЁхч $pr_x(a)$ ш $pr_y(a)$ яЁюхъЎшш Єюўъш
$a\in K$ эр юёш ъююЁфшэрЄ.

─ы  яюфьэюцхёЄтр $K\subset I^2$ юяЁхфхышь юЄюсЁрцхэшх
({\it ышэхщэ√щ юяхЁрЄюЁ ёєяхЁяючшЎшш})
$$\phi:C(I)\oplus C(I)\to
C(K)\quad\mbox{ЇюЁьєыющ}\quad\phi(g,h)(x,y)=g(x)+h(y).$$
═юЁьр эр яЁюёЄЁрэёЄтх $C(I)\oplus C(I)$ ттюфшЄё  хёЄхёЄтхээ√ь ёяюёюсюь.
▀ёэю, ўЄю яюфьэюцхёЄтю $K\subset I^2$ срчшёэю Єюуфр ш Єюы№ъю Єюуфр, ъюуфр
$\phi$ ¤яшьюЁЇэю.

╬сючэрўшь ўхЁхч $C^*(X)$ яЁюёЄЁрэёЄтю юуЁрэшўхээ√ї ышэхщэ√ї
ЇєэъЎшщ $C(X)\to\R$ ё эюЁьющ $|\mu|=\sup\{|\mu(f)|\ :\ f\in C(X),\ |f|=1\}$.
─ы  яюфьэюцхёЄтр $K\subset I^2$ юяЁхфхышь юЄюсЁрцхэшх
({\it фтющёЄтхээ√щ ышэхщэ√щ юяхЁрЄюЁ ёєяхЁяючшЎшш})
$$\phi^*:C^*(K)\to C^*(I)\oplus
C^*(I)\quad\mbox{ъръ}\quad \phi^*\mu(g,h)=(\mu(g\circ pr_x),\mu(h\circ pr_y)).$$
═юЁьр эр яЁюёЄЁрэёЄтх $C^*(I)\oplus C^*(I)$ ттюфшЄё  хёЄхёЄтхээ√ь ёяюёюсюь.
╥ръ ъръ $|\phi^*\mu|\le2|\mu|$, Єю $\phi^*$ юуЁрэшўхэ.
╧ю фтющёЄтхээюёЄш, $\phi$ ¤яшьюЁЇхэ Єюуфр ш Єюы№ъю Єюуфр, ъюуфр
$\phi^*$ ьюэюьюЁЇхэ.
\footnote{╧Ёш ¤Єюь $\phi^*$ ьюцхЄ с√Є№ шэ·хъЄштэ√ь, эю эх ьюэюьюЁЇэ√ь.
─Ёєушьш ёыютрьш, эх Єюы№ъю ышэхщэ√х ёююЄэю°хэш  эр $\im\phi$ чрёЄрты ■Є
хую с√Є№ ёЄЁюую ьхэ№°х ўхь $C(K)$, ъръ яюърч√трхЄ яЁшьхЁ эхсрчшёэющ
яюяюыэхээющ ьюыэшш.
\newline
╟рьхЄшь, ўЄю хёыш $K\subset\R^2$ - срчшёэюх яюфьэюцхёЄтю,
Єю ь√ ьюцхь фюърчрЄ№ схч шёяюы№чютрэш 
$\phi$, ўЄю $\phi^*$ ьюэюьюЁЇэю.
╬яЁхфхышь ышэхщэ√щ юяхЁрЄюЁ
$\Psi:C^*(I)\oplus C^*(I)\to C^*(K)$ ЇюЁьєыющ
$\Psi(\mu_x,\mu_y)(f)=\mu_x(g)+\mu_y(h)$,
уфх $g,h\in C(I)$ Єръют√, ўЄю $g(0)=0$ ш $f(x,y)=g(x)+h(y)$ фы 
$(x,y)\in K$. ▀ёэю, ўЄю
$\Psi\phi^*=\id$ ш $\Psi$ юуЁрэшўхэю, ёыхфютрЄхы№эю $\phi^*$ ьюэюьюЁЇэю.}

╧юэ Єэю, ўЄю $\phi^*$ ьюэюьюЁЇхэ Єюуфр ш
Єюы№ъю Єюуфр, ъюуфр

(*) {\it ёє∙хёЄтєхЄ  $\varepsilon>0$ Єръюх, ўЄю
$|\phi^*\mu|>\varepsilon|\mu|$ фы  ърцфюую эхэєыхтюую $\mu\in C^*(K)$.}

╫Єюс√ ЁрсюЄрЄ№ ё єёыютшхь (*), шёяюы№чєхь ёыхфє■∙шщ эхЄЁштшры№э√щ ЇръЄ:
{\it $C^*(K)$ ёютярфрхЄ ё яЁюёЄЁрэёЄтюь $\sigma$-aффшЄштэ√ї
Ёхуєы Ёэ√ї тх∙хёЄтхээючэрўэ√ї сюЁхыхтёъшї ьхЁ эр $K$}
(фрыхх ь√ сєфхь эрч√трЄ№ шї яЁюёЄю 'ьхЁрьш'; шёяюы№чєхЄё  Єръцх ЄхЁьшэ
'чрЁ ф√').
╚ьххь
$$\phi^*\mu=(\mu_x,\mu_y),\quad\mbox{уфх}\quad\mu_x(U)=\mu(pr_x^{-1}U)\quad
\mbox{ш}\quad\mu_y(U)=\mu(pr_y^{-1}U).$$
 ┼ёыш $\mu=\mu^+-\mu^-$ хёЄ№
Ёрчыюцхэшх ьхЁ√ $\mu$ эр яюыюцшЄхы№э√х ш юЄЁшЎрЄхы№э√х ўрёЄш, Єю
$|\mu|=\bar\mu(X)$, уфх $\bar\mu=\mu^++\mu^-$ хёЄ№ рсёюы■Єэюх чэрўхэшх
ьхЁ√ $\mu$.

─юърчрЄхы№ёЄтю Єюую, ўЄю єёыютшх (*) тыхўхЄ юЄёєЄёЄтшх ёъюы№ єуюфэю
фышээ√ї ьюыэшщ, юёЄрты хь т ърўхёЄтх єяЁрцэхэш .
(▌Єю фюърч√трхЄ ўрёЄ№ 'Єюы№ъю Єюуфр', фы  ъюЄюЁющ є эрё єцх хёЄ№ ¤ыхьхэЄрЁэюх
фюърчрЄхы№ёЄтю.)

╬ёЄр╕Єё  фюърчрЄ№, ўЄю єёыютшх (*) ёыхфєхЄ шч $E^n(K)=\emptyset$.
(▌Єю фюърч√трхЄ ўрёЄ№ 'Єюуфр'.)
╧Ёштхфхь фюърчрЄхы№ёЄтю фы  $n\in\{1,2\}$ (фы  яЁюшчтюы№эюую $n$ юэю
рэрыюушўэю).

╬сючэрўшь ўхЁхч $D_x$ (ш $D_y$) ьэюцхёЄтю Єхї Єюўхъ шч $K$, ъюЄюЁ√х эх
чрЄхэ ■Єё  эшъръющ фЁєующ Єюўъющ шч $K$ т $x$- (ш $y$-) эряЁртыхэшш.
┬юч№ьхь ы■сє■ ьхЁє $\mu$ эр $K$ ё эюЁьющ 1.

┼ёыш $n=1$, Єю
$$E(K)=\emptyset,\quad\mbox{яю¤Єюьє}\quad D_x\cup D_y=K,\quad\mbox{чэрўшЄ,}
\quad 1=\bar\mu(K)\le\bar\mu(D_x)+\bar\mu(D_y).$$
Tюуфр, эх єьхэ№°р  юс∙эюёЄш, $\bar\mu(D_x)\ge1/2$.
╥ръ ъръ $pr_x$ шэ·хъЄштэр эр $D_x$, Єю $|\mu_x|\ge1/2$.
╧ю¤Єюьє єёыютшх (*) т√яюыэхэю фы  $\varepsilon=\frac 12$.

┼ёыш $n=2$, Єю
$$E(E(K))=\emptyset,\quad\mbox{яю¤Єюьє}\quad D_x\cup D_y=K-E(K)\quad\mbox{ш}
\quad E(D_x\cup D_y)=\emptyset.$$
\quad
╧Ёш $\bar\mu(E(K))<3/4$ шьххь $\bar\mu(D_x\cup D_y)>1/4$ ш, эх єьхэ№°р 
юс∙эюёЄш, $\bar\mu(D_x)>1/8$, чэрўшЄ, ъръ ш т ёыєўрх $n=1$, шьххь
$|\mu_x|>1/8$.
╧ю¤Єюьє єёыютшх (*) т√яюыэхэю фы  $\varepsilon=\frac 18$.

╧Ёш $\bar\mu(E(K))\ge3/4$ шьххь $\bar\mu(K-E(K))\le1/4$.
╩ръ ш т ёыєўрх $n=1$, эх єьхэ№°р  юс∙эюёЄш,
$\bar\mu_x(pr_x(E(K)))\ge\bar\mu(E(K))/2$.
╤ыхфютрЄхы№эю $|\mu_x|\ge\frac12\cdot\frac34-\frac14=\frac18$.
╧ю¤Єюьє єёыютшх (*) т√яюыэхэю фы  $\varepsilon=\frac 18$.
QED

\normalsize

\bigskip
{\bf ├ырфър  срчшёэюёЄ№.}

╧єёЄ№ $K$ --- яюфьэюцхёЄтю яыюёъюёЄш $\R^2$.
╘єэъЎш  $f:K\to\R$ эрч√трхЄё  {\it фшЇЇхЁхэЎшЁєхьющ} (яю ╙шЄэш), хёыш фы  ы■сющ Єюўъш
$z_0\in K$ ёє∙хёЄтє■Є Єръшх тхъЄюЁ $a\in\R^2$ ш схёъюэхўэю ьрыр  ЇєэъЎш 
$\alpha:\R^2\to\R$, ўЄю фы  ы■сющ Єюўъш $z\in K$ т√яюыэхэю
$$f(z)=f(z_0)+a\cdot(z-z_0)+\alpha(z-z_0)|z,z_0|.$$
\quad
╟фхё№ Єюўър ючэрўрхЄ чэръ ёъры Ёэюую яЁюшчтхфхэш  тхъЄюЁют $a=:(f_x,f_y)$ ш
$z-z_0=:(x,y)$, Є.х. $a\cdot(z-z_0)=xf_x + yf_y$.
╘єэъЎш  $\alpha:\R^2\to\R$
эрч√трхЄё  {\it схёъюэхўэю ьрыющ}, хёыш фы  ы■сюую ўшёыр $\eps >0$ ёє∙хёЄтєхЄ
Єръюх ўшёыю $\delta>0$, ўЄю фы  ы■сющ Єюўъш $(x,y)\in \R^2$
$$\mbox{хёыш}\quad \sqrt{x^2+y^2}<\delta,\quad\mbox{Єю}\quad|\alpha(x,y)|<\eps.$$

╧юфьэюцхёЄтю $K\subset\R^2$ яыюёъюёЄш эрч√трхЄё  {\it фшЇЇхЁхэЎшЁєхью срчшёэ√ь},
хёыш фы  ы■сющ фшЇЇхЁхэЎшЁєхьющ ЇєэъЎшш $f:K\to\R$ ёє∙хёЄтє■Є Єръшх
фшЇЇхЁхэЎшЁєхь√х ЇєэъЎшш $g:\R\to\R$ ш $h:\R\to\R$, ўЄю фы  ы■сющ Єюўъш
$(x,y)\in K$ т√яюыэ хЄё  $f(x,y)=g(x)+h(y)$.

\smallskip
{\bf 13.} (a) (b) (c) ╨х°шЄх рэрыюуш чрфрўш 6 фы  фшЇЇхЁхэЎшЁєхьющ срчшёэюёЄш.

\smallskip
{\bf 14.}
(a) ├ЁрЇшъ ЇєэъЎшш $|x|$ эр юЄЁхчъх $[-1;1]$  ты хЄё  фшЇЇхЁхэЎшЁєхью срчшёэ√ь.

(b) ╦юьрэр  ё яюёыхфютрЄхы№э√ьш тхЁ°шэрьш $(-2,0)$, $(-1,1)$, $(0,0)$, $(1,1)$
ш $(2,0)$ эх  ты хЄё  фшЇЇхЁхэЎшЁєхью срчшёэющ.
(╟рьхЄшь, ўЄю юэр срчшёэр.)


(c) ╧юяюыэхээр  ьюыэш 
$\{([\frac{n+1}2]^{-1/2},[\frac n2]^{-1/2})\}_{n=1}^{\infty}\cup\{(0,0)\}$
эх  ты хЄё  фшЇЇхЁхэЎшЁєхью срчшёэющ.
(╟рьхЄшь, ўЄю юэр эх  ты хЄё  Єръцх срчшёэющ.)

(d) ╧юяюыэхээр  ьюыэш 
 $\{(2^{-[\frac{n+1}2]},2^{-[\frac n2]})\}_{n=1}^{\infty}\cup\{(0,0)\}$
 ты хЄё  фшЇЇхЁхэЎшЁєхью срчшёэющ.
(╟рьхЄшь, ўЄю юэр эх  ты хЄё  срчшёэющ.)

(e)** {\bf ├шяюЄхчр ╚. ╪эєЁэшъютр.} 
╧юяюыэхээр  ьюыэш  $\{a_n\}_{n=1}^\infty\cup\{(0,0)\}$ 
фшЇЇхЁхэЎшЁєхью срчшёэр Єюуфр ш Єюы№ъю Єюуфр, ъюуфр яюёыхфютрЄхы№эюёЄ№ 
$\frac{\sum\limits_{n=k}^\infty|a_n|}{a_k}$ юуЁрэшўхэр. 

\smallskip
{\bf 15.}
(a) ╩ЁхёЄ $K=[(-1,-2),(1,2)]\cup[(-1,1),(1,-1)]$ эх  ты хЄё  фшЇЇхЁхэЎшЁєхью
срчшёэ√ь.

(b)** {\bf ├шяюЄхчр.} ╧юфьэюцхёЄтю
$\{(t^2,\frac {t^2}{(1+t)^2})\}_{t\in[-\frac 12;\frac 12]}$ яыюёъюёЄш
эх  ты хЄё  фшЇЇхЁхэЎшЁєхью срчшёэ√ь.
╙ърчрэшх: ьюцэю я√ЄрЄ№ё  фхырЄ№ рэрыюушўэю чрфрўх 15(a).

(c)** {\bf  ├шяюЄхчр.}
╩єёюўэю-ышэхщэ√щ уЁрЇ эр яыюёъюёЄш  ты хЄё  фшЇЇхЁхэЎшЁєхью срчшёэ√ь Єюуфр ш 
Єюы№ъю Єюуфр, ъюуфр юэ эх ёюфхЁцшЄ ёъюы№ єуюфэю фышээ√ї ьюыэшщ, ш фы  ы■с√ї 
фтєї {\it ёшэуєы Ёэ√ї} Єюўхъ $a$ ш $b$ т√яюыэхэю $x(a)\ne x(b)$ ш 
$y(a)\ne y(b)$.
╥юўър $a\in K$ эрч√трхЄё  {\it ёшэуєы Ёэющ}, хёыш яхЁхёхўхэшх $K$
ё ы■с√ь фшёъюь ё ЎхэЄЁюь т $a$ эх  ты хЄё  яЁ ьюышэхщэ√ь юЄЁхчъюь.

(d)** ═рщфшЄх ъЁшЄхЁшщ фшЇЇхЁхэЎшЁєхьющ срчшёэюcЄш фы  уЁрЇют т яыюёъюёЄш.

(e)** ╤є∙хёЄтєхЄ ыш эхяЁхЁ√тэюх юЄюсЁрцхэшх юЄЁхчър т яыюёъюёЄ№,
юсЁрч ъюЄюЁюую  ты хЄё  фшЇЇхЁхэЎшЁєхью срчшёэ√ь, эю эх срчшёэ√ь?

\smallskip
{\bf 16.}
╧єёЄ№ $K$ --- яюфьэюцхёЄтю яыюёъюёЄш $\R^2$ ш $r\ge0$.
╘єэъЎш  $f:K\to \R$ эрч√трхЄё  {\it $r$ Ёрч фшЇЇхЁхэЎшЁєхьющ}, хёыш фы  ы■сющ
Єюўъш $z_0\in K$ ёє∙хёЄтє■Є Єръшх ьэюуюўыхэ $\overline f(z)=\overline f(x,y)$
ёЄхяхэш эх т√°х $r$ юЄ фтєї яхЁхьхээ√ї $x$ ш $y$ ш схёъюэхўэю ьрыр  ЇєэъЎш 
$\alpha:\R^2\to\R$, ўЄю
$f(z)=\overline f(z-z_0)+\alpha(z-z_0)|z,z_0|^r$ фы  ы■сющ Єюўъш $z\in K$.
(▌Єю юяЁхфхыхэшх юЄышўрхЄё  юЄ юс∙хяЁшэ Єюую.)

(a) ═юы№ Ёрч фшЇЇхЁхэЎшЁєхь√х ЇєэъЎшш --- ¤Єю т ЄюўэюёЄш эхяЁхЁ√тэ√х, р
юфшэ Ёрч фшЇЇхЁхэЎшЁєхь√х --- ¤Єю т ЄюўэюёЄш фшЇЇхЁхэЎшЁєхь√х.

(b) ─ы  ы■сюую Ўхыюую яюыюцшЄхы№эюую $r$ юяЁхфхышЄх $r$-фшЇЇхЁхэЎшЁєхьє■
 срчшёэюёЄ№ яюфьэюцхёЄт яыюёъюёЄш.

(c)* ─ы  ы■сюую Ўхыюую $k\ge0$ эрщфeЄё  яюфьэюцхёЄтю яыюёъюёЄш,
$r$-фшЇЇхЁхэЎшЁєхью срчшёэюх фы  ы■сюую $r=0,1,\dots k$, эю эх
$r$-фшЇЇхЁхэЎшЁєхью срчшёэюх эш фы  ъръюую $r>k$ [RZ07].

(d)** ═рщфшЄх ъЁшЄхЁшщ $r$-фшЇЇхЁхэЎшЁєхьющ срчшёэюcЄш фы  уЁрЇют т яыюёъюёЄш.

\bigskip
{\bf ╨х°хэш  чрфрў.}

\small
\smallskip
{\bf 13.}
(a), (b), (c) └эрыюушўэю чрфрўрь 6(a), 3(a) ш 3(b).

\smallskip
{\bf 14.}
(a) ╧єёЄ№ $f(x,y)$ - фшЇЇхЁхэЎшЁєхьр  ЇєэъЎш . ╥юуфр
$f(x,|x|)-f(0,0)=ax + b|x| +\alpha(x,|x|)|(x,|x|)|$.
╧юыюцшь $h(y)=by$, $g(x)=f(x,|x|)-b|x|$.
╧юфЁюсэхх ёь. [RZ07].

(b) ╧Ёхфяюыюцшь, ўЄю фрээр  ыюьрэр  фшЇЇхЁхэЎшЁєхью срчшёэр.
╘єэъЎш  $f(x,y)=xy$ фшЇЇхЁхэЎшЁєхьр. 
╧ю¤Єюьє $f(x,y)=g(x)+h(y)$ фы  эхъюЄюЁ√ї фшЇЇхЁхэЎшЁєхь√ї ЇєэъЎшщ $g$ ш $h$.
╥юуфр 
$$2-2d=f(1+d,1-d)+f(1-d,1-d)=g(1+d)+g(1-d)+2h(1-d)=
2g(1)+2h(1)-2h'(1)d+o(d).$$
╟эрўшЄ, $h'(1)=1$. 
└эрыюушўэю  
$$2d-2=f(-1+d,1-d)+f(-1-d,1-d)=g(-1+d)+g(-1-d)+2h(1-d)=
2g(-1)+2h(1)-2h'(1)d+o(d).$$
╟эрўшЄ, $h'(1)=-1$. 
╧ЁюЄштюЁхўшх. 

(c) ╧Ёхфяюыюцшь, ўЄю ¤Єр яюяюыэхээр  ьюыэш  фшЇЇхЁхэЎшЁєхью срчшёэр. ╧юыюцшь
\linebreak
 $a_n=([\frac{n+1}2]^{-1/2},[\frac n2]^{-1/2})$,
$f(a_n):=\frac{(-1)^n}n$,
$n=2,3,\dots$.
┼ёыш $f(x,y)=g(x)+h(y)$ фы  эхъюЄюЁ√ї ЇєэъЎшщ $g(x)$ ш $h(y)$,
Єюуфр
$f(a_2)-f(a_3)+f(a_4)-\dots$ ёїюфшЄё  ъ $g(1)-g(0)$
(рэрыюушўэю чрфрўх 7b). ═ю ¤Єю яЁюЄштюЁхўшЄ ЁрёїюфшьюёЄш Ё фр
$\frac 12 +\frac 13+ \frac 14 + \dots $

(d) ═х юуЁрэшўштр  юс∙эюёЄш ьюцэю ёўшЄрЄ№, ўЄю $f(0,0)=0$, Єюуфр тюч№ь╕ь
$g(0)=0$ ш $h(0)=0$.
╧юыюцшь
$$h(2^{-k})=f(2^{-(k+1)},2^{-k})-f(2^{-(k+1)},2^{-(k+1)})+f(2^{-(k+2)},2^{-(k+1)})-
\dots, $$
$$g(2^{-k})=f(2^{-k},2^{-k})-f(2^{-(k+1)},2^{-k})+
f(2^{-(k+1)},2^{-(k+1)})-\dots,$$
уфх яЁрт√х ўрёЄш ёєЄ№ ёєьь√ чэръюяхЁхьхээ√ї Ё фют.

╥хяхЁ№ $g(x)$ ш $h(y)$ ьюуєЄ с√Є№ яЁюфюыцхэ√ фю фшЇЇхЁхэЎшЁєхь√ї ЇєэъЎшщ
 $\R\to\R$.

\smallskip
{\bf 15.} (a)
╬яЁхфхышь
$$w(0)=0,\quad w(4^{-i}+4^{-3i})=w(4^{-i})=0\quad\mbox{ ш }
\quad w(4^{-i}+4^{-3i-1})=2^{3i}\quad\mbox{ фы  }\quad i=1,2,3,\dots.$$
╧Ёюфюыцшь ЄхяхЁ№ ¤Єє ЇєэъЎш■ ъєёюўэю-ышэхщэю фю ЇєэъЎшш $w:[0;1]\to \R$.
─ы  $x\in[0;1]$ юяЁхфхышь $W(x)$ ъръ яыю∙рф№ яюф уЁрЇшъюь
ЇєэъЎшш $w$ эр юЄЁхчъх $[0;x]$.
╬яЁхфхышь $f(x,-x)=W(x)$ фы  $x\in[0;1]$ and $f(x,y)=0$ эр юёЄры№э√ї Єюўърї 
ъЁхёЄр. 

▀ёэю, ўЄю $f$ фшЇЇхЁхэЎшЁєхьр тэх $(0,0)$. 
╠юцэю яЁютхЁшЄ№, ўЄю $f$ фшЇЇхЁхэЎшЁєхьр ш т $(0,0)$. 

╧Ёхфяюыюцшь, ўЄю $f(x,y)=g(x)+h(y)$ фы  эхъюЄюЁ√ї фшЇЇхЁхэЎшЁєхь√ї
$g$ ш $h$. 
═х юуЁрэшўштр  юс∙эюёЄш, сєфхь ёўшЄрЄ№, ўЄю $g(0)=h(0)=0$. 
╘єэъЎш  $g$ эх фшЇЇхЁхэЎшЁєхьр т Єюўъх $x=1/4$, яюёъюы№ъє фы  $0<d<\frac 14$ 
т√яюыэхэю 
$$g\left(\frac14+d\right)-g\left(\frac14\right)=
W\left(\frac14+d\right)-W\left(\frac14\right)+
W\left(\frac1{4^2}+\frac d4\right)-W\left(\frac1{4^2}\right)+\dots >$$
$$>W\left(\frac1{4^{k+1}}+\frac d{4^k}\right)-W\left(\frac1{4^{k+1}}\right)=
\frac{2^{3k}\cdot 4^{-3k}}2\ge \frac {(4d)^{3/4}}2.$$
╟фхё№  

$\bullet$ яхЁтюх ЁртхэёЄтю фюърч√трхЄё  ё шёяюы№чютрэшхь фтєї схёъюэхўэ√ї 
ьюыэшщ шч Єюўхъ ъЁхёЄр, эрўшэр■∙шїё  т Єюўърї $(\frac14+d,-\frac14-d)$
ш $(\frac14,-\frac14)$ ш ёїюф ∙шїё  ъ Єюўъх $(0,0)$; 

$\bullet$ $k\ge0$ Єръютю, ўЄю $4^{-2k}\ge 4d>4^{-2(k+1)}$;  

$\bullet$ яхЁтюх эхЁртхэёЄтю т√яюыэхэю, яюёъюы№ъє $W$ эхтючЁрёЄр■∙р  ЇєэъЎш ; 

$\bullet$ тЄюЁюх эхЁртхэёЄтю т√яюыэхэю, яюёъюы№ъє 
$\frac d{4^k}>\frac1{4^{3(k+1)}}$; 

$\bullet$ тЄюЁюх ЁртхэёЄтю т√яюыэхэю яю юяЁхфхыхэш■ ўшёыр $k$.

(└эрыюушўэю фюърч√трхЄё , ўЄю $g$ эх фшЇЇхЁхэЎшЁєхьр эш т ъръющ Єюўъх тшфр
$x=4^{-i}$.)

\normalsize



\bigskip
\centerline{\uppercase{\bf ╦шЄхЁрЄєЁр}}

[Ar58] ┬. ╚. └Ёэюы№ф,
{\it ╬ яЁхфёЄртыхэшш ЇєэъЎшщ эхёъюы№ъшї яхЁхьхээ√ї т тшфх ёєяхЁяючшЎшш ЇєэъЎшщ
ьхэ№°хую ўшёыр яхЁхьхээ√ї,} ╠рЄ. ╧Ёюётх∙хэшх, ╤хЁ.~2, т√я.~3, (1958), 41--61.
http://ilib.mirror1.mccme.ru/djvu/mp2/mp2-3.djvu?, djvuopts\&page=43

[Ar58'] ┬. ╚. └Ёэюы№ф, {\it ╧Ёюсыхьр 6,} ╠рЄ. ╧Ёюётх∙хэшх, ╤хЁ.~2,
т√я.~3, (1958), 273-274.
\linebreak
http://ilib.mirror1.mccme.ru/djvu/mp2/mp2-3.djvu?, djvuopts\&page=243

[Vi04] └. ├. ┬шЄє°ъшэ. \emph{13-  яЁюсыхьр ├шы№схЁЄр ш ёьхцэ√х тюяЁюё√,}
╙╠═, Є.~59, т√я.~1(355), 2004. ╤.~11--24.

[Vo81] ╤. ╠. ┬юЁюэшэ. \emph{└эрышЄшўхёър  ъырёёшЇшърЎш  ЁюёЄъют ъюэЇюЁьэ√ї
юЄюсЁрцхэшщ $(\mathbb C,0)\to(\mathbb C,0)$ ё ЄюцфхёЄтхээющ ышэхщэющ
ўрёЄ№■,} ╘єэъЎ. рэрышч ш хую яЁшы. ╥.~15, т√я.~1, 1981. ╤.~1--17.

[Vo82] ╤. ╠. ┬юЁюэшэ. \emph{└эрышЄшўхёър  ъырёёшЇшърЎш  ярЁ шэтюы■Ўшщ ш хх
яЁшыюцхэш ,} ╘єэъЎ. рэрышч ш хую яЁшы. ╥.~16, т√я.~2, 1982. ╤.~21--29.

[Ku68] ╩. ╩єЁрЄютёъшщ, ╥юяюыюуш , ╠шЁ, ╠юёътр, 1969, Є. 1,2

[Ku00] V.~Kurlin, {\it Basic embeddings into products of graphs,}
Topol.Appl. 102 (2000), 113--137.

[Ku03] ┬. └. ╩єЁышэ. \emph{┴рчшёэ√х тыюцхэш  уЁрЇют ш ьхЄюф ЄЁхїёЄЁрэшўэ√ї
тыюцхэшщ ─√ээшъютр,} ╙╠═, Є.~58, т√я.~2(350), 2003. ╤.~163--164.

[Ku03'] ┬. └. ╩єЁышэ,  ┴рчшёэ√х тыюцхэш  уЁрЇют ш ьхЄюф ЄЁхїёЄЁрэшўэ√ї
тыюцхэшщ ─√ээшъютр, фшёёхЁЄрЎш  (2003),
http://maths.dur.ac.uk/$\sim$dma0vk/PhD.html

[Mi09] E. Miliczka, {\it Constructive decomposition of a function of two
variables as a sum of functions of one variable}, Proc. AMS, 137:2 (2009),
607-614.

[MKT03] N. Mramor-Kosta and E. Trenklerova [Miliczka],
{\it On basic embeddings of
compacta into the plane,} Bull. Austral. Math. Soc. 68 (2003), 471--480.


[Pr04] ┬. ┬. ╧Ёрёюыют, ▌ыхьхэЄ√ ъюьсшэрЄюЁэющ ш фшЇЇхЁхэЎшры№эющ Єюяюыюушш,
╠юёътр, ╠╓═╠╬, 2004. http://www.mccme.ru/prasolov


[RS99]
─. ╨хяют° ш └. ╤ъюяхэъют, ''═ют√х Ёхчєы№ЄрЄ√ ю тыюцхэш ї яюыш¤фЁют ш
ьэюуююсЁрчшщ т хтъышфют√ яЁюёЄЁрэёЄтр'', ╙╠═ 54:6 (1999), p. 61--109

[RZ07] D. Repovs and M. Zeljko, {\it On basic embeddings into the plane,}
Rocky Mountain J. Math., to appear.

[Sk95] A. Skopenkov, {\it A description of continua basically embeddable in
$\R^2$,} Topol. Appl. 65 (1995), 29--48.

[Sk05] └. ╤ъюяхэъют, ┬юъЁєу ъЁшЄхЁш  ╩єЁрЄютёъюую яырэрЁэюёЄш уЁрЇют,
╠рЄ. ╧Ёюётх∙хэшх, 9 (2005), 116--128 ш 10 (2006), 276--277.
http://www.mccme.ru/free-books/matprosa.html,
http://arxiv.org/abs/0802.3820

[Sk] └. ┴. ╤ъюяхэъют, └ыухсЁршўхёър  Єюяюыюуш  ё ¤ыхьхэЄрЁэющ Єюўъш чЁхэш ,
╠юёътр, ╠╓═╠╬, т яхўрЄш.
http://arxiv.org/abs/0808.1395

[St89] Y.~Sternfeld, {\it Hilbert's 13th problem and dimension,}
Lect. Notes Math. 1376 (1989), 1--49.

\end{document}